\documentclass{amsart}

\usepackage{epsf,latexsym,amsfonts,amssymb,subfigure,mathrsfs,graphicx,stmaryrd}
\usepackage{longtable}
\usepackage{color}
\newcommand{\secref}[1]{\S~\ref{#1}}
\newcommand{\abs}[1]{\lvert#1\rvert}
\newcommand{\labs}[1]{\left\lvert\,#1\,\right\rvert}

\newcommand{\Lr}[1]{\left(#1\right)}
\newcommand{\lr}[1]{\bigl(#1\bigr)}
\newcommand{\llr}[1]{\left[#1\right]}
\newcommand{\set}[2]{\{\,#1\,\mid\,#2\,\}}
\newcommand{\diff}[2]{\dfrac{\partial #1}{\partial #2}}

\newcommand{\ffd}[1]{D_{#1}^+}

\newcommand{\mc}[1]{\mathcal{#1}}
\newcommand{\mb}[1]{\mathbb{#1}}

\def\lot{\text{L.O.T.}}
\def\Om{\Omega}

\def\x{\times}

\def\al{\alpha}
\def\del{\delta}
\def\ga{\gamma}
\def\eps{\varepsilon}
\def\vr{\varrho}

\def\na{\nabla}

\def\lam{\lambda}

\def\F{\boldsymbol F}

\def\0{\boldsymbol 0}
\def\1{\boldsymbol 1}
\def\Z{\mathbb Z}
\def\R{\mathbb R}

\def\md{\mathrm{d}}

\def\dx{\md\,x}

\def\co{\mathcal O}

\def\F{\mathcal F}

\def\Q{\mathcal Q}

\def\wcb{W_{\text{CB}}}

\def\yycb{y_{\text{cb}}}
\def\yat{y_{\text{at}}}
\def\fat{\F_{\text{at}}}

\def\fqc{\F_{\text{qc}}}

\def\etot{E^{\text{tot}}}

\def\wcb{W_{\text{cb}}}

\newcommand{\nn}{\nonumber}

\def\argmin{\operatorname{argmin}}
\def\argmax{\operatorname{argmax}}

\newtheorem{theorem}{Theorem}[section]
\newtheorem{lemma}[theorem]{Lemma}
\newtheorem{corollary}[theorem]{Corollary}

\newtheorem{remark}[theorem]{Remark}

\theoremstyle{definition}

\numberwithin{equation}{section}
\begin {document}
\title[QC Analysis]
{Ghost Force Influence of a Quasicontinuum Method in Two Dimension\vrule height
15pt width 0pt}
\author[J. Chen]{Jingrun Chen}
\address{Institute of
Computational Mathematics and Scientific/Engineering Computing,
AMSS, Chinese Academy of Sciences, No. 55, Zhong-Guan-Cun East Road,
Beijing, 100190, China}\email{chenjr@lsec.cc.ac.cn}
\curraddr{South Hall 6705, Mathematics Department,
University of California, Santa Barbara, CA93106, USA}\email{cjr@math.ucsb.edu}

\author[P.-B. Ming]{Pingbing Ming}
\address{LSEC, Institute of
Computational Mathematics and Scientific/Engineering Computing,
AMSS, Chinese Academy of Sciences, No. 55 Zhong-Guan-Cun East Road,
Beijing, 100190, China} \email{mpb@lsec.cc.ac.cn}

\thanks{The work of Ming was partially supported by National Natural Science
              Foundation of China grants 10932011, by the funds from
              Creative Research Groups of China through grant 11021101, and by the
              support of CAS National Center for Mathematics and Interdisciplinary
              Sciences. The authors are grateful to Weinan E and Jianfeng Lu for helpful discussions.}

\date{\today}
\begin {abstract}
We derive an analytical expression for the solution of a two-dimensional quasicontinuum method with a planar interface. The
expression is used to prove
that the ghost force may lead to a finite size error for the gradient of the solution. We estimate
the width of the interfacial layer
induced by the ghost force is of $\co(\sqrt{\eps}\,)$ with $\eps$ the
equilibrium bond length, which is much wider than that of the one-dimensional
problem.
\end {abstract}

\subjclass[2000]{65N30, 65N12, 65N06, 74G20, 74G15}
\keywords{Quasicontinuum method, ghost force, interfacial layer}
\maketitle
\section{Introduction}

Multiscale methods have been developed to simulate mechanical behaviors of
solids for several
decades~\cite{MillerTadmor:2009}. Combination of models at different
scales greatly enhances the dimension of problems that computers can deal with.
However, problems regarding the consistency,
stability and convergence of the multiscale methods may arise from the
coupling~\cite {E:book}. Taking the quasicontinuum (QC)
method~\cite{TadmorOrtizPhillips:1996, KnapOrtiz:2001} for example, one of the main issues is the so called ghost force
problem~\cite{ShenoyMillerTadmorPhillipsOrtiz:1999}, which is the artificial
non-zero force that the atoms experience at the
equilibrium state. In the language of numerical analysis,
the scheme lacks consistency at the interface between the atomistic region
and the continuum region~\cite{ELuYang:2006}.  For the one-dimensional problem,
it has been shown in~\cite{MingYang:2009} that the ghost force may lead to
a finite size error for the gradient of the solution. This generates
an interfacial layer with width $\co(\eps\abs{\ln\eps})$, out of which the error
for the gradient of the solution is of $\co(\eps)$.

To understand the influence of the ghost force for high dimensional problems, we study
a two-dimensional triangular lattice model with a QC approximation. This
QC method couples the Cauchy-Born
elasticity model and the atomistic model with a planar interface.
Numerical results show that
the ghost force may lead to a finite size error
for the gradient of the solution as in the
one-dimensional problem. The error profile
exhibits a layer-like structure. It seems that
the width of the interfacial layer induced by the ghost force is of
$\co(\sqrt{\eps})$, which is much wider than that of the one-dimensional problem.

To further characterize the influence of the ghost force, we introduce a
square lattice model with a QC approximation. Compared to
the triangular lattice model, this model can be solved {\em analytically} and
the error profile exhibits a clear layer-like structure.
Based on the analytical solution, we prove the error committed by the ghost force
for the gradient of the solution is $\co(1)$ and
the width of the induced interfacial layer is of $\co(\sqrt{\eps})$, which are
also confirmed by the numerical results.

The paper is organized as follows. Numerical results for
the triangular model and the square lattice model
with QC approximations are presented in~\secref{sec:tri}
and~\secref{sec:square}, respectively.
We derive an analytical expression of the solution of the square lattice model with
a QC approximation in~\secref{sec:explicit}.
The main results of
the paper are proved in~\secref{sec:estimate}.
%
\section{A quasicontinuum method for triangular lattice}~\label {sec:tri}
\subsection{Atomistic and continuum models}
\begin{figure}
\centering
\includegraphics[width=2.5in]{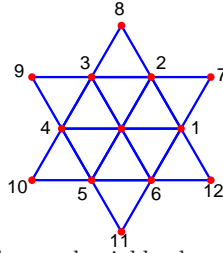}%
\setlength{\abovecaptionskip}{-30pt}\caption{\small The first and
second neighborhood interaction ranges of the triangular lattice; $\mc{S}_1=\{s_1,\cdots,s_6\}$
and $\mc{S}_2=\{s_7,\cdots,s_{12}\}$.}\label {stencilschematic}
\end{figure}

We consider the triangular lattice $\mb{L}$, which can be written as
\[
\mb{L}=\set{x\in\R^2}{x=ma_1+na_2,\;m,n\in\Z}
\]
with the basis vectors
$a_1=(1,0),a_2=(1/2,\sqrt3/2)$.
Define the unit cell of $\mb{L}$ as
\[
\Gamma=\set{x\in\R^2}{x=c_1a_1+c_2a_2,\;-1/2\le c_1,c_2<1/2}.
\]
We shall consider lattice system $\eps\mb{L}$ inside the domain $\Om$,
and $\Om_\eps=\Om\cap\eps\mb{L}$, where $\eps$ is the equilibrium
bond length. Assume
that the atoms are interacted with the potential function $V$,
which is usually a highly nonlinear function, e.g., the
Lennard-Jones potential~\cite {LennardJones:1924}.
Denote by $\mc{S}_1$ and $\mc{S}_2$ the first and the second neighborhood
interaction ranges; see Figure~\ref{stencilschematic}. In particular, we have
\begin{align*}
\mc{S}_1&=\cup_{i=1}^6s_i=\{a_1,a_2,-a_1+a_2,-a_1,-a_2,a_1-a_2\},\\
\mc{S}_2&=\cup_{i=7}^{12}s_i
=\{a_1+a_2,-a_1+2a_2,-2a_1+a_2,-a_1-a_2,a_1-2a_2,2a_1-a_2\}.
\end{align*}

For $\mu\in\Z^2$, the translation
operator $T^{\mu}_{\eps}$ is defined
for any lattice function $z:\mb{L}\to\R^2$ as
\[
  (T^{\mu}_{\eps}z)(x) = z(x + \eps \mu_1a_1+\eps\mu_2a_2)\quad \text{for } x \in \mb{L}.
\]
We define the forward and backward discrete gradient operators as
\[
D_s^+ =\eps^{-1}(T^{\mu}_{\eps}-I)\qquad\text{and}\qquad
D_s^-=\eps^{-1}(I-T^{\mu}_{\eps}),
\]
where $s =\mu_1a_1+\mu_2a_2$ and $I$ is the identity operator. We shall
also use the short-hand
$Dz=(D_1^+z, D_2^+z)=(D_{s_1}^+z,D_{s_2}^+ z)$.

Consider an atomic system posed on $\Om_\eps$. The total energy is given by
\begin{equation}\label{eq:atomener}
\etot_{\text{at}}=\dfrac12\sum_{x\in\Om_\eps}\sum_{s\in\mc{S}_1\cup\mc{S}_2}
V\Lr{\abs{\ffd{s}z(x)}},
\end{equation}
where $V$ is a potential function. In this paper, we only consider the pairwise potential function, and leave the discussion on the more general potential
functions in future publication.
The atomistic problem is to minimize the total energy subject to certain
boundary conditions that will be specified later on.

Next we turn to the Cauchy-Born elasticity
model~\cite {BornHuang:1954, Ericksen:1984, EMing:2007, Ericksen:2008}. Given a $2$ by $2$ matrix $A$, the stored energy density function
is given by
\[
\wcb(A)=\dfrac{1}{2\vartheta_0}\sum_{s\in\mc{S}_1\cup\mc{S}_2}V(\abs{s\cdot
A}),
\]
where $\vartheta_0$ is the area of the unit cell and
$\vartheta_0=\sqrt{3}\eps^2/2$. The stored energy function is defined by
\[
\etot_{\text{cb}}=\int_\Om\wcb(\na z(x))\,\dx.
\]
The continuum problem is to minimize the stored energy function
subject to certain boundary conditions. We employ the standard $P_1$
Lagrange finite element to approximate the Cauchy-Born elasticity model
with the lattice $\mb{L}$ as the
triangulation. The approximate stored energy function is
\begin{equation}\label {eq:v-ener}
\etot_{\text{cb},\eps}=\dfrac12\sum_{x\in\Om_\eps}\sum_{i=1}^6\Lr{
V\Lr{\abs{\ffd{s_i}z(x)}}+V\Lr{\abs{(\ffd{s_i}+\ffd{s_{i+1}})z(x)}}}.
\end{equation}
One can see $\etot_{\text{cb},\eps}$ reproduces the atomistic energy $\etot_{\text{at}}$; cf.~\eqref{eq:atomener}, if only the nearest neighborhood interaction is considered.

\begin{figure}
\centering
\includegraphics[width=3in]{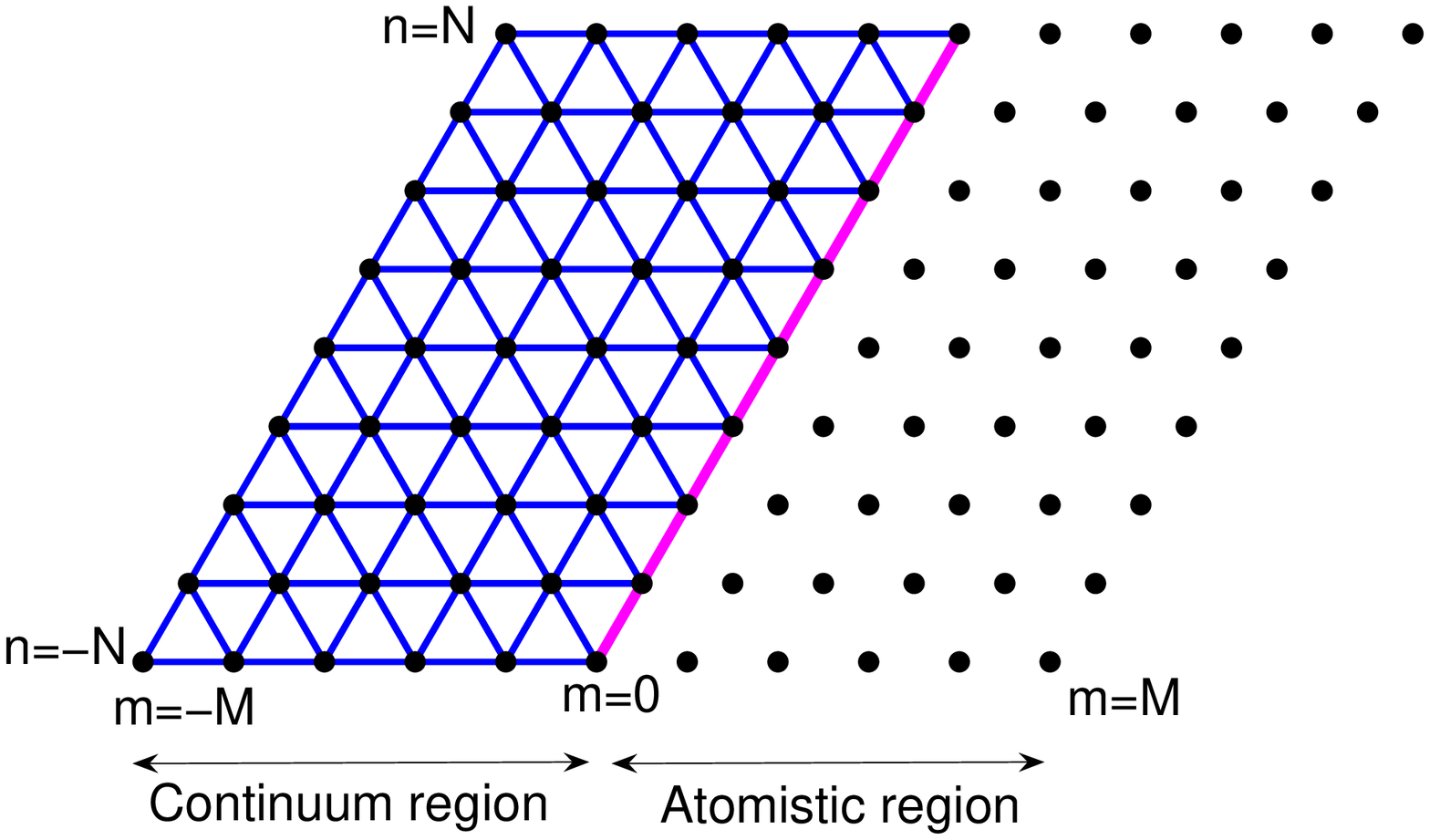}
\setlength{\abovecaptionskip}{-15pt}
\caption{\small Schematic picture of $\Om_\eps$.}\label {schematic}
\end{figure}

We study the quasicontinuum method~\cite {TadmorOrtizPhillips:1996}.
Let $\eps=1/(2M)$. We assume that the
interface between the continuum model and the atomistic model
is $x_1=0$ as shown in Figure~\ref{schematic}.
The total energy of the QC method is
\begin{align*}
\etot_{\text{qc}} & =\dfrac12\sum_{x_1\leq -2\eps}\sum_{i=1}^6\Lr{
V\Lr{\abs{\ffd{s_i}z(x)}}+V\Lr{\abs{(\ffd{s_i}+\ffd{s_{i+1}})z(x)}}}\\
&\quad +\dfrac12\sum_{x_1= -\eps} \bigg\{ \sum_{i=1}^6\Lr{
V\Lr{\abs{\ffd{s_i}z(x)}}+V\Lr{\abs{(\ffd{s_i}+\ffd{s_{i+1}})z(x)}}}
+ V(\abs{\ffd{s_{12}}z(x)}) \bigg\}\\
&\quad +\dfrac12\sum_{x_1= 0} \bigg\{\sum_{s\in\mc{S}_1}
V\Lr{\abs{\ffd{s}z(x)}} + \sum_{i=2}^4\Lr{
V\Lr{\abs{\ffd{s_i}z(x)}}+V\Lr{\abs{(\ffd{s_i}+\ffd{s_{i+1}})z(x)}}}\\
&\phantom{+\dfrac12\sum_{x_1= 0} \bigg\{}
\qquad +\dfrac12V(\abs{\ffd{s_{7}}z(x)})
+ V(\abs{\ffd{s_{11}}z(x)})+ V(\abs{\ffd{s_{12}}z(x)}) \bigg\}\\
&\quad +\dfrac12\sum_{x_1= \eps} \bigg\{ \sum_{s\in\mc{S}_1\cup\mc{S}_2}
V\Lr{\abs{\ffd{s}z(x)}} - \dfrac12V(\abs{\ffd{s_{9}}z(x)}) \bigg\}\nn \\
&\quad +\dfrac12\sum_{x_1\ge 2\eps}\sum_{s\in\mc{S}_1\cup\mc{S}_2}
V\Lr{\abs{\ffd{s}z(x)}}.
\end{align*}
The force at atom $x$ is defined by
\[
\fqc[z](x)\equiv-\diff{\etot_{\text{qc}}}{z(x)}.
\]

Since we only concern the influence of the ghost force,
following~\cite {MingYang:2009}, we
assume that the potential function is a harmonic function, i.e.,
\[
V(r)=\dfrac12{r^2}.
\]
Without taking into account the external force, we write
the equilibrium equations for the QC approximation as
\[
\fqc[z](x)=0
\]
with
\begin{align*}
\fqc[z](x)&=-12 z(x)+\sum_{i=1}^{12}z(x+\eps s_i),&& x\in\Om_\eps, x_1\le -2\eps,\\
\fqc[z](x)&=-24 z(x)+4\sum_{i=1}^6z(x+\eps s_i),&&x\in\Om_\eps, x_1\ge 2\eps.
\end{align*}

For $x=(-\eps,x_2)$,
\[
\fqc[z](x)=-\dfrac{49}{2}z(x)+4\sum_{i=1}^6z(x+\eps s_i)
+\dfrac12 z(x+\eps s_{12}).
\]

For $x=(0,x_2)$,
\begin{align*}
\fqc[z](x)&=-18z(x)+[z(x+\eps s_1)+z(x+\eps s_6)]+\dfrac52[z(x+\eps s_2)
+z(x+\eps s_5)]\\
&\quad+4[z(x+\eps s_3)+z(x+\eps s_4)]+[z(x+\eps s_7)+z(x+\eps s_{11})+z(x+\eps s_{12})].
\end{align*}

For $x=(\eps,x_2)$,
\[
\fqc[z](x)=-\dfrac{23}{2}z(x)+\sum_{i=1}^{12}z(x+\eps s_i)-\dfrac12
z(x+\eps s_{9}).
\]

At the equilibrium state, we evaluate $\fqc$ at $z(x)=x$ to get
\[
\fqc[x](x)=\left\{
\begin{aligned}
(-3\eps/4,\sqrt3\eps/4),&\rm{\quad if\quad}x_1=-\eps,\\
(3\eps/2,-\sqrt{3}\eps/2),&\rm{\quad if\quad}x_1=0,\\
(-3\eps/4,\sqrt3\eps/4),&\rm{\quad if\quad}x_1=\eps,\\
(0,0),&\rm{\quad otherwise}.
\end{aligned}\right.
\]
The above equations imply $z_1(x)=-\sqrt3 z_2(x)$. Therefore,
we only study the first component of $z(x)$ and neglect the subscript
if no confusion will occur. The error of the displacement induced by
the ghost force, denoted by $y(x)\equiv z(x)-x$, satisfies
\begin{equation}\label{eq:qctrierr}
\fqc[y](x)=\fqc[z](x)-\fqc[x](x)=-\fqc[x](x)\equiv f(x)
\end{equation}
with
\[
f(x)=\left\{
\begin{aligned}
3\eps/4,&\text{\quad if\quad}x_1=-\eps,\\
-3\eps/2,&\rm{\quad if\quad}x_1=0,\\
3\eps/4,&\rm{\quad if\quad}x_1=\eps,\\
0&\rm{\quad otherwise}.
\end{aligned}\right.
\]

Boundary conditions need to be supplemented to close the system of equilibrium
equations. Two types of boundary conditions will be considered in this paper.
One is the periodic boundary condition in $x_2$ direction
while homogeneous Dirichlet boundary condition in $x_1$ direction, which will be called
periodic boundary condition if no confusion will occur. The
other is the homogeneous Dirichlet boundary conditions in both $x_1$ and
$x_2$ directions.

In what follows, we shall use a conventional notation
$y(m,n)=y(x)$ with $x=\eps(ma_1+na_2)$.
\subsection{Periodic boundary condition}
The periodic boundary condition can be written as
\[
\left\{
\begin{aligned}
y(-M,n)&=0,&&n=-N,\cdots,N,\\
y(m,n)&=y(m,2N+n),&&m=-M,\cdots,M, n=-N,\cdots,N,\\
y(M,n)&=y(M+1,n)=0,&&n=-N,\cdots,N.
\end{aligned}\right.
\]

Observe that the solution of~\eqref{eq:qctrierr} with the periodic
boundary condition takes a special form:
\[
y(x)=cy(x_1),
\]
namely, the solution is constant along $x_2$ direction. Based on this observation, we
conclude that the QC approximation with the periodic boundary condition
reduces to a one-dimensional problem. The
equilibrium equations satisfied by $y(x_1)$ are as follows.
In the continuum region, i.e., $m=-M+1,\cdots,-2$,
\[
8y(m+1)-16y(m)+8y(m-1)=0,
\]
and in the atomistic region, i.e., $m=2,\cdots,M-1$,
\[
y(m+2)+4y(m+1)-10y(m)+4y(m-1)+y(m-2)=0.
\]
The equations in the interfacial region are
\[
\left\{
\begin{aligned}
\dfrac12 y(1)+8 y(0)-\dfrac{33}{2} y(-1)+8y(-2)&=\dfrac34\eps,\\
y(2)+4 y(1)-13 y(0)+8 y(-1)&=-\dfrac32\eps,\\
y(3)+4 y(2)-\dfrac{19}{2} y(1)+4 y(0)+\dfrac12y(-1)&=\dfrac34\eps.
\end{aligned}\right.
\]
The boundary condition is
\[
y(-M)=y(M)=y(M+1)=0.
\]

The above equilibrium equations can also be obtained by considering a
one-dimensional chain interacted with the following harmonic potential:
\begin{equation}\label{harmonicpotential}
V(\{y\}) = \dfrac{k_1}{2}\sum_{\abs{i-j}=1}\abs{y_i-y_j}^2
+\dfrac{k_2}{2}\sum_{\abs{i-j}=2}\abs{y_i-y_j}^2
\end{equation}
with $k_1=4$ and $k_2=1$. This is the model studied
in~\cite{DobsonLuskin:2009a}.
\subsection{Dirichlet boundary condition}\label {tri_dbc}
\begin{figure}[htbp]
\centering
\subfigure[$D_{s_1}^+ y$]{%
\includegraphics[width=2.5in]{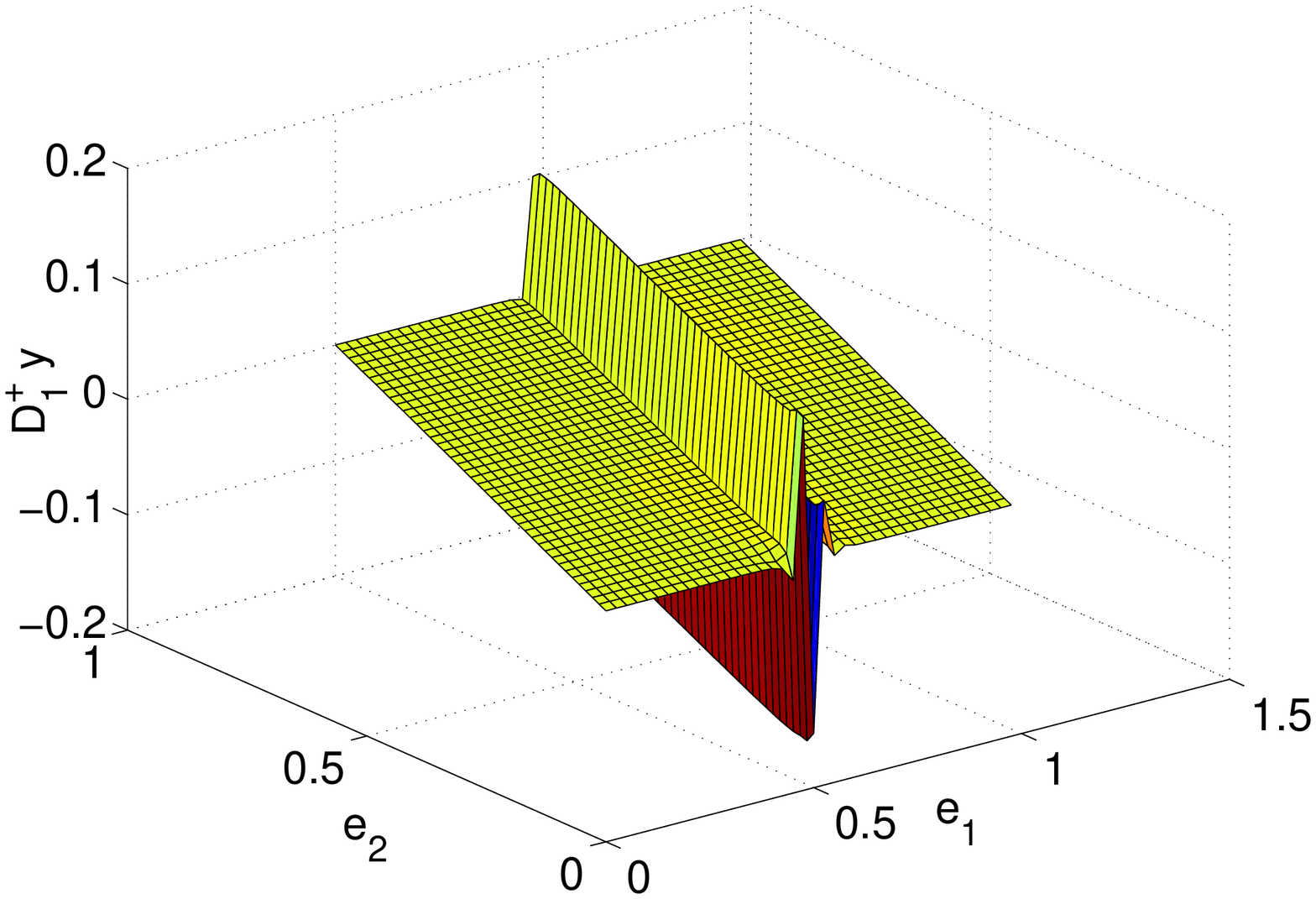}}%
\subfigure[$D_{s_2}^+ y$]{
\includegraphics[width=2.5in]{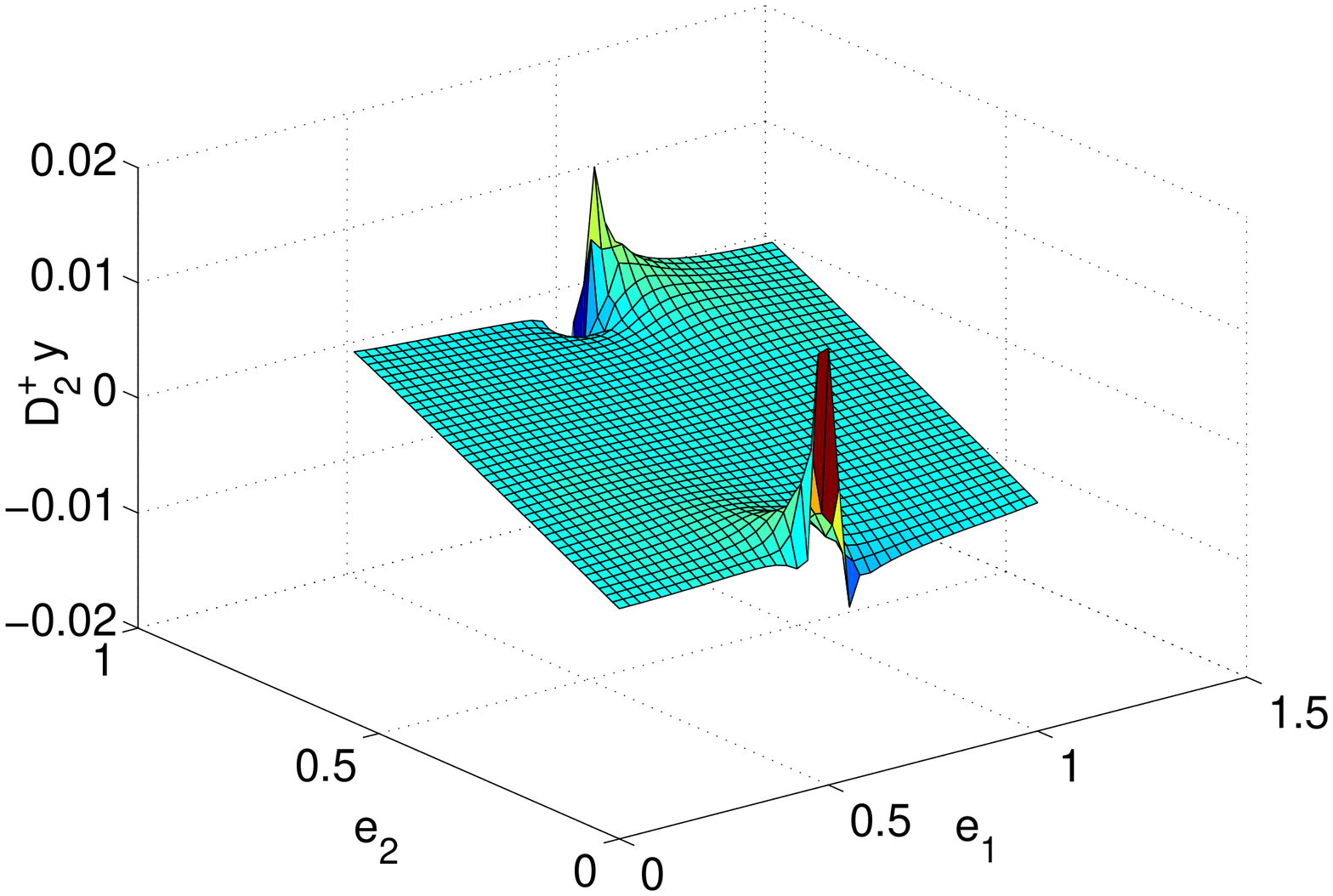}}%
\caption{\small Profiles of the discrete gradients for triangular lattice with $M=N=20$ under
Dirichlet boundary condition. }\label {triDBCinterface}
\end{figure}

\begin{figure}[htbp]
\centering
\subfigure[$\chi_{E_1}(x)$]{%
\label {triinterface1}
\includegraphics[width=2.5in]{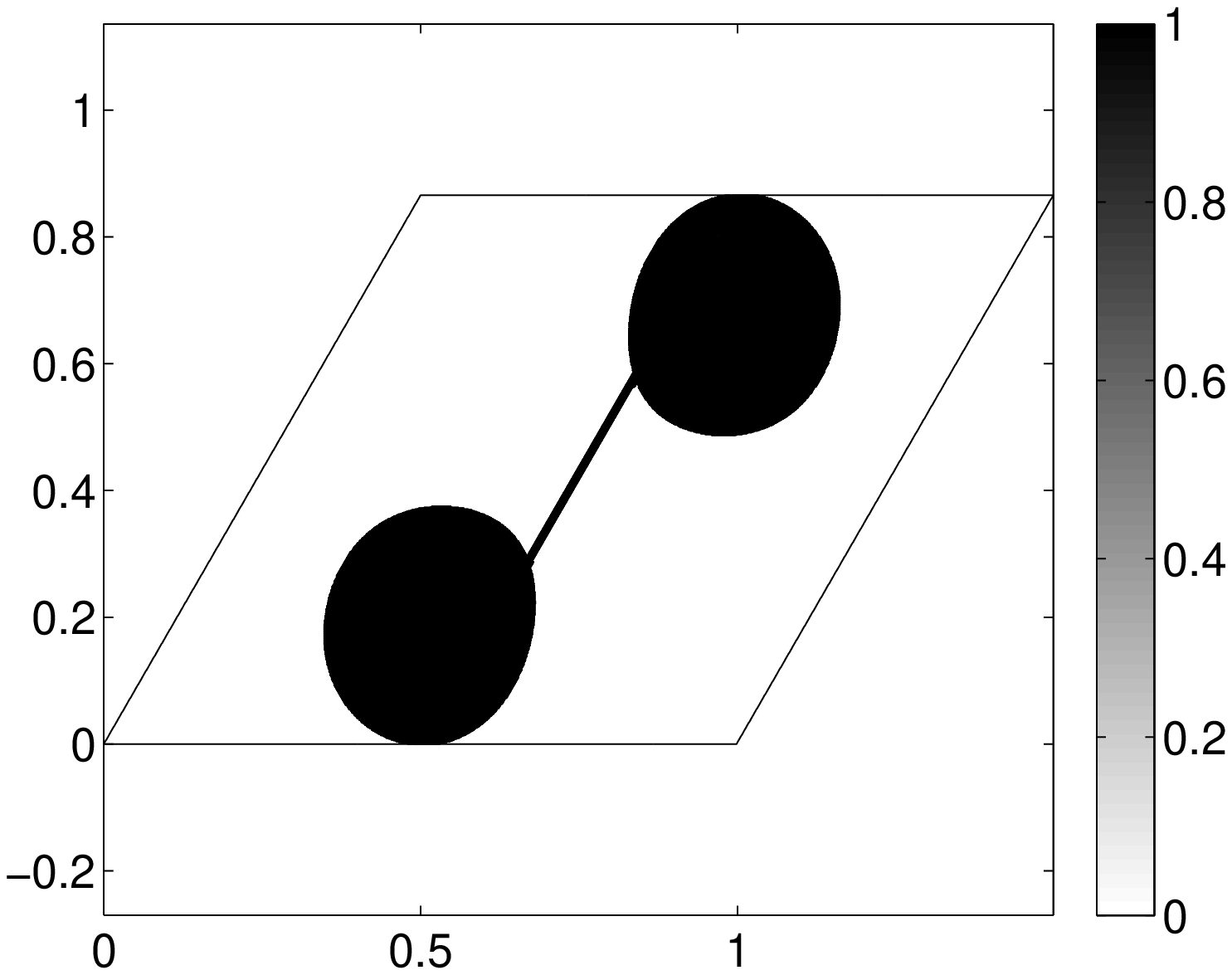}}%
\subfigure[$\chi_{E_2}(x)$]{
\label {triinterface2}
\includegraphics[width=2.5in]{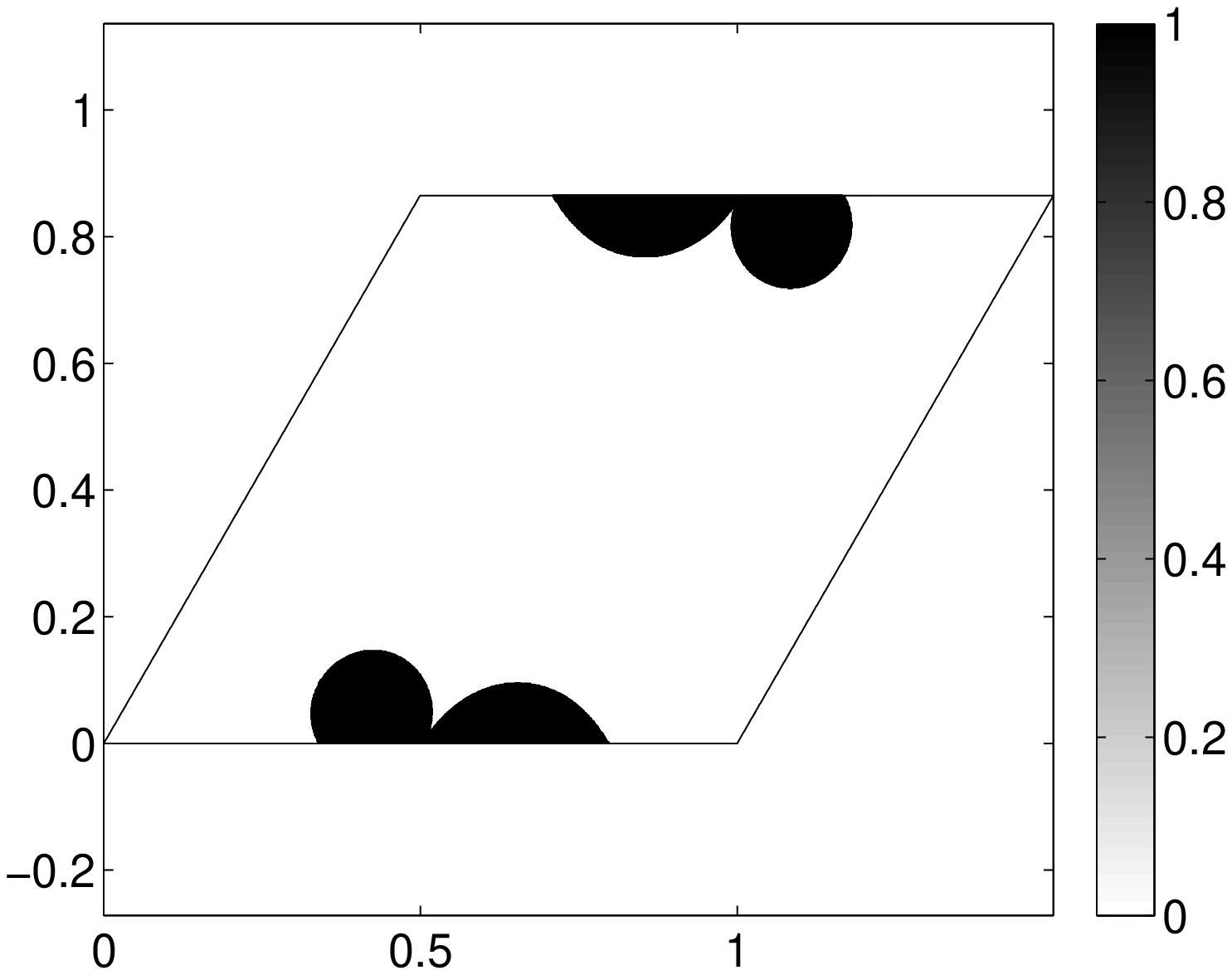}}%
\caption{\small  Profiles of the characteristic functions for triangular lattice.}\label {triinterface}
\end{figure}

The Dirichlet boundary condition can be written as
\[
\left\{
\begin{aligned}
y(-M,n)&=y(M,n)=y(M+1,n)=0,&& n=-N,\cdots,N,\\
y(m,-N)&=y(m,N)=0,&& m=-M,\cdots,-1,\\
y(m,-N-1)&=y(m,-N)=y(m,N)=y(m,N+1)=0,&& m=0,\cdots,M.
\end{aligned}\right.
\]
This gives an essentially two-dimensional model. We show profiles of discrete gradients
$D_{s_1}^+y$ and $D_{s_2}^+y$ in
Figure~\ref{triDBCinterface} with $M=N=20$. The profile of $D_{s_1}^+y$
is similar to that of the one-dimensional problem. However, if one zooms in the interface,
there are some differences.
To highlight the interface,  for $i=1,2$, we define
\[
E_i=\set{x\in\Om_\eps}{\abs{D_{s_i}y(x)}\ge c_0\eps},
\]
and let $\chi_{E_i}$ be their characteristic functions.
We plot $\chi_{E_i}(x)$ in Figure~\ref{triinterface} with $M=N=640$ and
$c_0=0.04$. In Figure~\ref{triinterface1}, the width of the
interface near the boundary is much wider than that in the interior domain.
No interfacial layer for $D_{s_2}^+y$ is observed in Figure~\ref{triinterface2}.
Therefore, we only measure the
width of the interfacial layer for $D_{s_1}^+ y$, which is defined by
\begin{align*}
\max_{n} \Big\{|a_1|\eps & \big[\argmax_{m}\{x=ma_1\eps+na_2\eps\in\Om_\eps| \chi_1(x) =1\} \\
& -\argmin_{m}\{x=ma_1\eps+na_2\eps\in\Om_\eps| \chi_1(x) =1\}\big]\Big\}.
\end{align*}
Numerical results in Table~\ref {trilayer} imply that the width of the interfacial layer
scales $\co(\sqrt{\eps})$, while the one-dimensional problem has an interfacial layer with
$\co(\eps\abs{\ln\eps})$ width~\cite {MingYang:2009,DobsonLuskin:2009a}.
\begin{table}[htbp]
\centering
\begin{tabular} {cccccccc}\hline
$\eps\;(5\x10^{-2})$  & $2^0$  & $2^{-1}$ & $2^{-2}$ & $2^{-3}$ & $2^{-4}$ & $2^{-5}$ & \text{Rate}\\
\hline \text{Layer width}$(10^{-1})$  & $3.5$ & $2.3$ & $1.5$ & $0.94$ & $0.63$ & $0.52$ & $0.57$\\
\hline
\end{tabular}
\vskip .3cm
\caption{\small Width of the interfacial layer versus the equilibrium bond
length $\eps$ for the triangular lattice with Dirichlet boundary condition.} \label {trilayer}
\end{table}

\begin{figure}
\centering
\includegraphics[width=3in,height=2in]{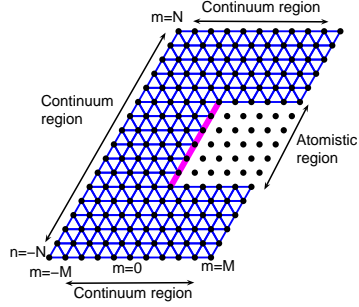}
\setlength{\abovecaptionskip}{-15pt}
\caption{{\small Padding technique to remove the boundary effect.}}\label {triRemove}
\end{figure}

One may doubt the above result could be caused by the boundary condition instead of the
ghost force. To clarify this issue, we enlarge the continuum region to
weaken the influence of the boundary condition.
We set $N=3M$ and use different equilibrium equations for different regions
as in Figure~\ref{triRemove}. Roughly speaking, along $x_2$ direction,
the original QC method is employed in the interior region, and the
continuum model is padded in the outer region.

\begin{figure}[htbp]
\centering
\subfigure[$D_{s_1}^+ y$]{%
\includegraphics[width=2.5in]{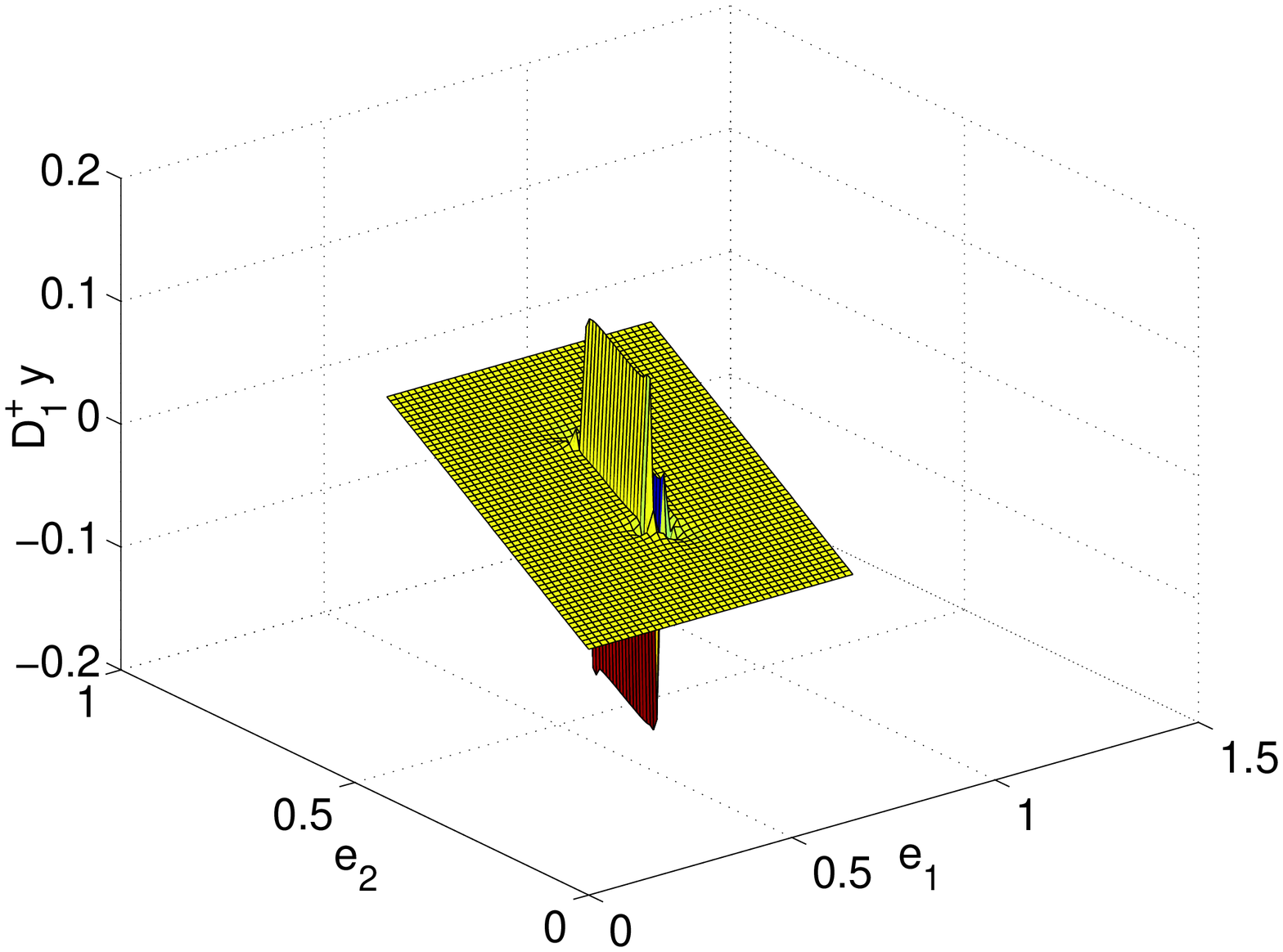}}%
\subfigure[$D_{s_2}^+ y$]{
\includegraphics[width=2.5in]{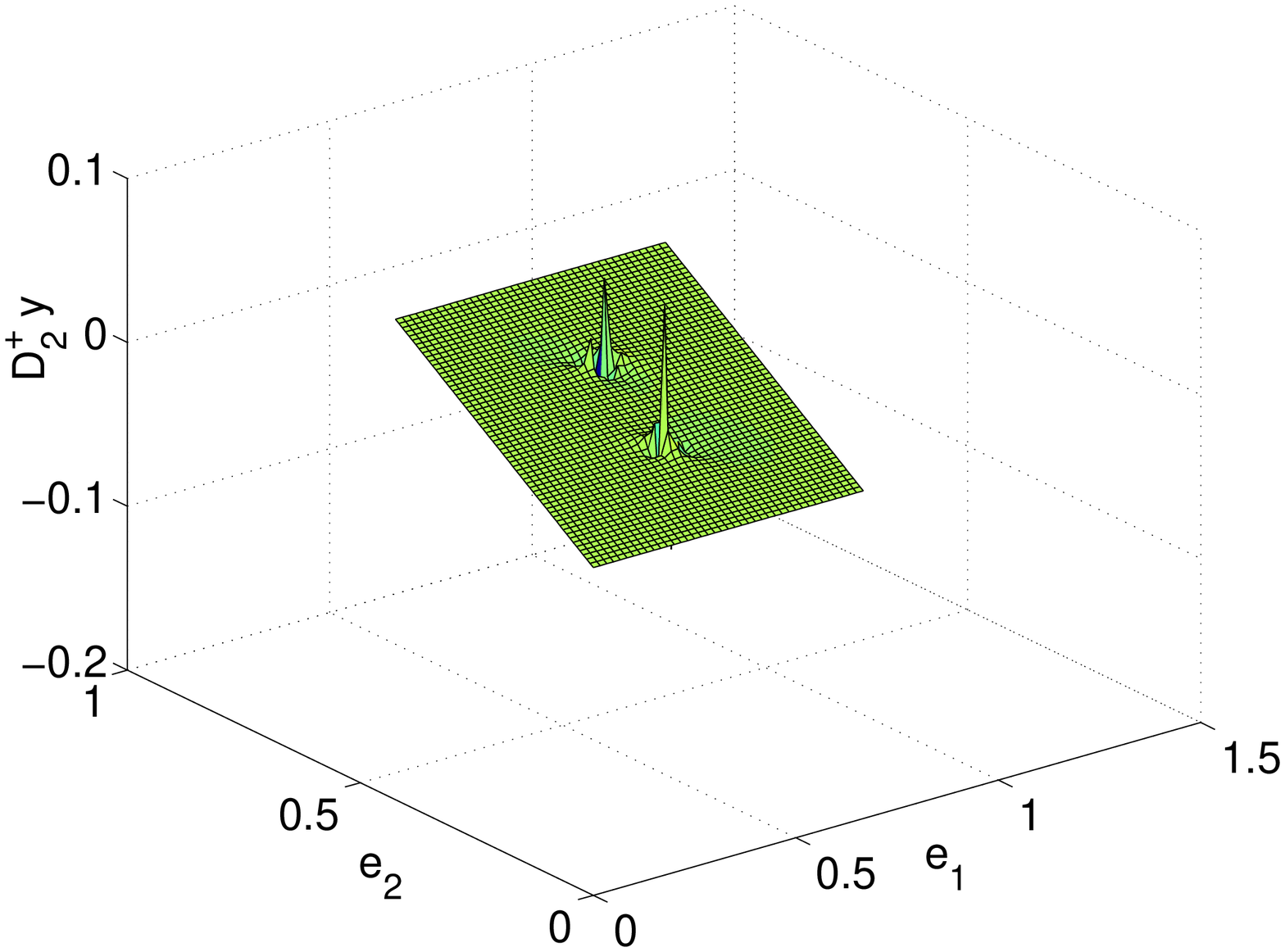}}%
\caption{\small Profiles of the discrete gradients for triangular lattice with $M=20$ and padding
technique under Dirichlet boundary condition. }\label {triDBCinterfaceCB}
\end{figure}

\begin{figure}[htbp]
\centering
\subfigure[$\chi_{E_1}(x)$]{%
\label {triinterfaceCB1}
\includegraphics[width=2.5in]{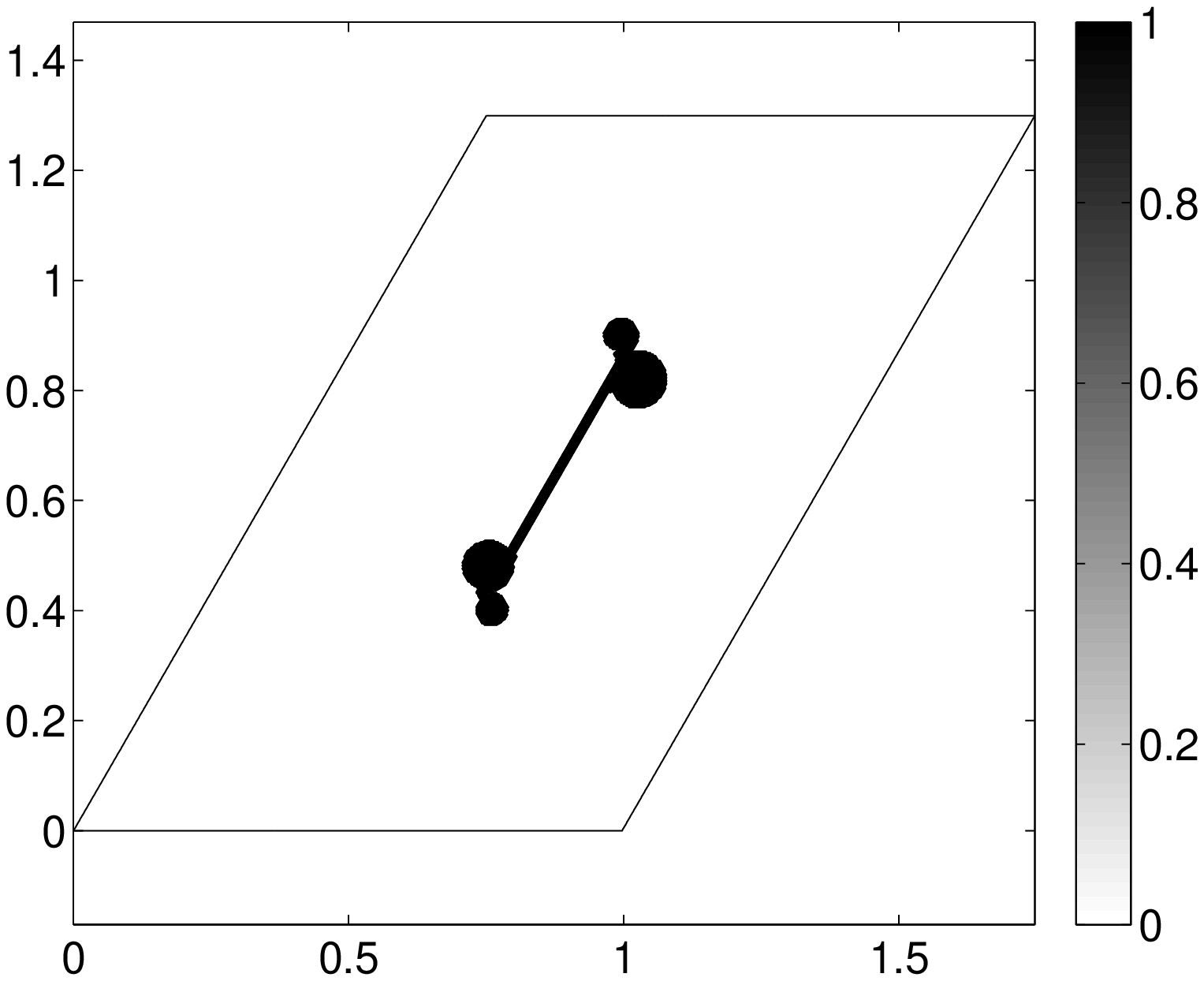}}%
\subfigure[$\chi_{E_2}(x)$]{
\label {triinterfaceCB2}
\includegraphics[width=2.2in]{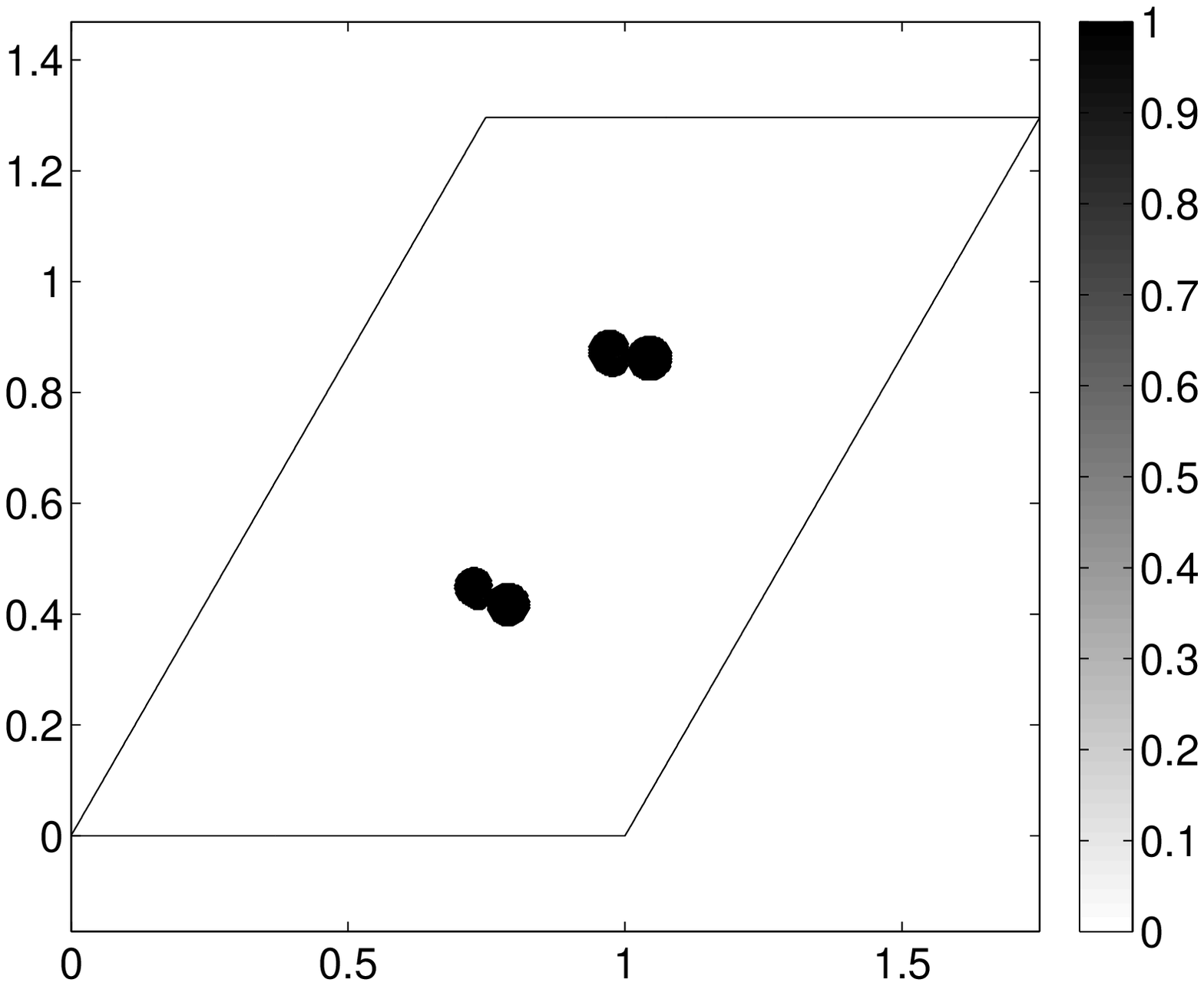}}%
\caption{\small  Profiles of the characteristic functions for triangular lattice
with padding technique.}\label {triinterfaceCB}
\end{figure}

We plot $D^+_{s_1} y$, $D^+_{s_2} y$, and their characteristic functions
in Figure~\ref{triDBCinterfaceCB} and Figure~\ref{triinterfaceCB}. They
are similar to those in Figure~\ref{triDBCinterface} and Figure~\ref{triinterface}.

Table~\ref {triRlayerCB}
shows the width of the interfacial layer in terms of $\eps$.
Numerical results suggest that the width of the interfacial
layer is still of $\co(\sqrt{\eps}\,)$, which is consistent with that in
Table~\ref{trilayer}. Therefore, we conclude that the interface is
caused by the ghost force instead of the boundary condition.
\begin{table}[htbp]
\centering
\begin{tabular} {cccccccc}\hline
$\eps(10^{-2})$  & $1.7$  & $0.83$ & $0.41$ & $0.21$ & $0.10$  & \text{Rate}\\
\hline \text{Layer width}$(10^{-2})$  & $9.2$ & $8.3$ & $4.8$ & $3.3$ & $2.3$ & $0.53$\\
\hline
\end{tabular}
\vskip .3cm
\caption{\small Width of the interfacial layer versus the equilibrium
bond length $\eps$ by
removing the effect of the boundary condition.}\label {triRlayerCB}
\end {table}

The QC approximation discussed in this section is quite realistic except
the potential function as shown in~\cite {MingYang:2009}, which
however seems enough to characterize the influence of the ghost force. Unfortunately, 
it does not seem easy to solve this model analytically as we have done 
in~\cite {MingYang:2009} for the one-dimensional problem. In next section, we 
introduce a new QC method
that can be solved {\em analytically}. We shall prove that this new QC method 
does capture the main feature of the ghost force for the triangular lattice model 
with a planar interface.
\section{Square lattice model}~\label {sec:square}
\begin{figure}
\centering
\begin{minipage}[c]{0.49\textwidth}
\centering
\rotatebox{0}{\scalebox{.85}{\includegraphics[width=2.8in]{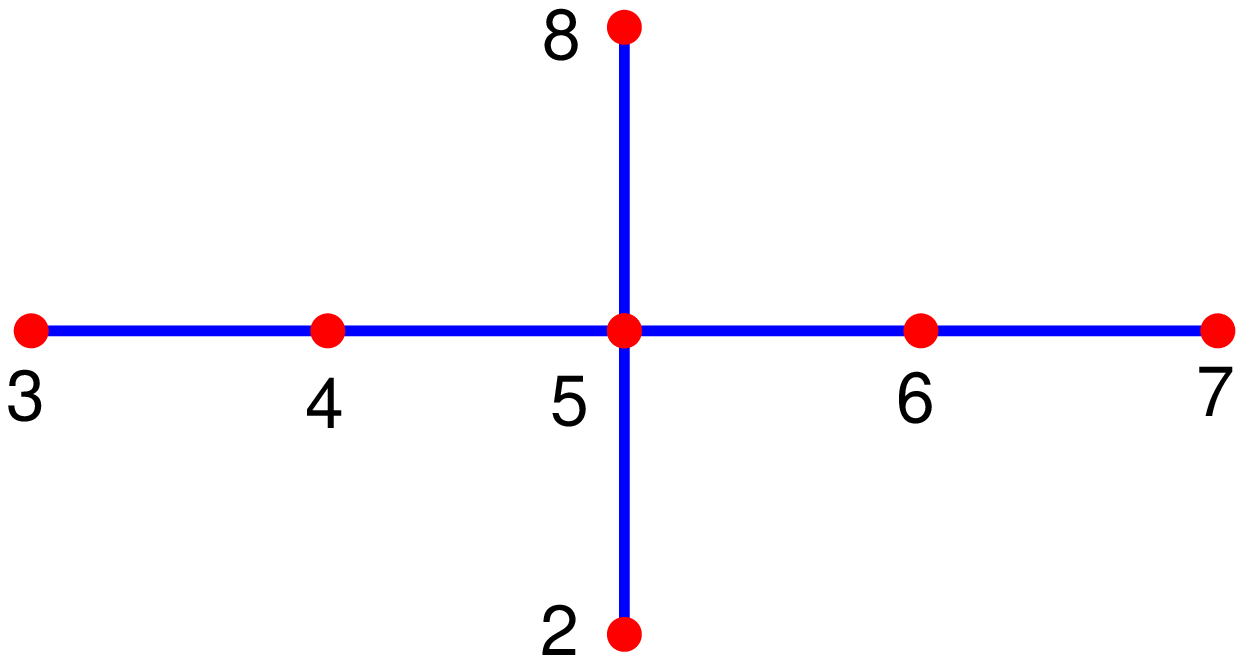}}}
\end{minipage}
\begin{minipage}[c]{0.49\textwidth}
\centering
\rotatebox{0}{\scalebox{.85}{\includegraphics[width=2.8in]{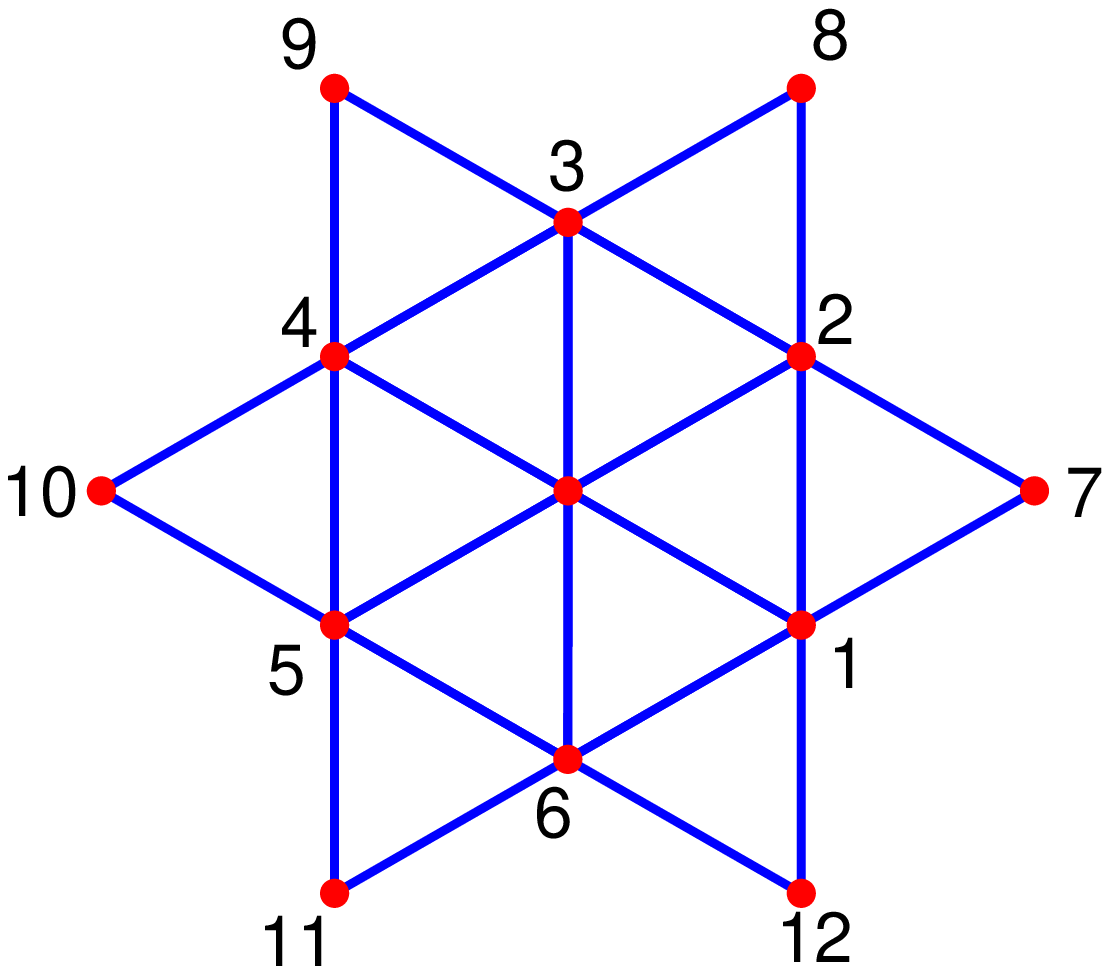}}}
\end{minipage}
\caption{\small Interaction ranges. Left: First and second neighborhood
interactions in $x_1$ direction and first neighborhood interaction in
$x_2$ direction for square lattice; Right: First and second neighborhood
interactions for triangular lattice in a rotated coordinate.}\label {stencil}
\end{figure}

We consider the square lattice with the
harmonic potential. Compared to the standard interaction range of
the square lattice, we assume a special interaction range
as shown in Figure~\ref{stencil} (Left). Namely, the
first and second neighborhood interactions in $x_1$ direction, and
the first neighborhood interaction in $x_2$ direction are taken into account.
This seemingly strange selection may be obtained from a rotated triangular lattice as in
Figure~\ref {stencil} (Right). If we condense the
interaction of the atoms $9,4,5, 11$
into one atom, and the atoms $8,2,1,12$ into another, then we obtain
a square lattice model with the special interaction range described above,
which may be the underlying reason why it can be regarded as a {\em surrogate} model.
We shall show in the next two sections that
this model not only captures the main features of the triangular lattice model with the QC approximation as shown in Figure~\ref{schematic}, but also lends
itself theoretically tractable.

Proceeding along the same line that leads to~\eqref {eq:qctrierr}, we obtain
the equilibrium equations for the error $y(x)$.\footnote{We actually multiply $-1$
on both sides of~\eqref{eq:qctrierr}.} In the continuum
region, i.e., $m=-M,\cdots,-2$ and $n=-N,\cdots,N$,
\begin{align}
12y_{i}(m,n) &
-y_{i}(m,n-1)-y_{i}(m,n+1)-5y_{i}(m-1,n)-5y_{i}(m+1,n)\nn\\
&= 0, \quad i=1,2,\label {eq:cb}
\end{align}
and in the atomistic region,  i.e., $m=2,\cdots,M$ and $n=-N,\cdots,N$,
\begin{align}
6y_{i}(m,n) &
-y_{i}(m,n-1)-y_{i}(m,n+1)-y_{i}(m-1,n)-y_{i}(m+1,n)\nn\\
&-y_{i}(m-2,n)-y_{i}(m+2,n)= 0,\quad i=1,2.\label {eq:atom}
\end{align}

\begin{figure}
\centering
\includegraphics[width=3in]{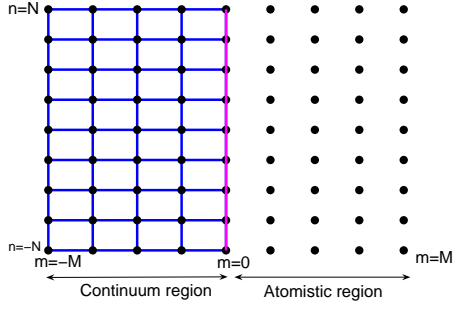}
\caption{\small Square lattice with a planar interface.}\label {planar}
\end{figure}

The interface between the continuum model and the atomistic model is the line $m=0$ as shown in
Figure~\ref {planar}, and $M$ is assumed to be even for simplicity.
The equilibrium equations for the layers $m=-1, 0$ and $1$ are as follows.

For layer $m=-1$ and $n=-N,\cdots,N$,
\begin{align}
\frac{25}2y_{i}(-1,n) &
-y_{i}(-1,n-1)-y_{i}(-1,n+1)-5y_{i}(-2,n)-5y_{i}(0,n)\nn\\
&-\frac12y_{i}(1,n)= f_{i}, \quad i=1,2,\label{eq:inters1}
\end{align}
where $f_{1}=-\eps$ and $f_{2}=0$.

For layer $m=0$ and $n=-N,\cdots,N$,
\begin{align}
9y_{i}(0,n) &
-y_{i}(0,n-1)-y_{i}(0,n+1)-5y_{i}(-1,n)-y_{i}(1,n)\nn\\
&-y_{i}(2,n)= f_{i}, \quad i=1,2,\label{eq:inters2}
\end{align}
where $f_{1}=2\eps$ and $f_{2}=0$.

For layer $m=1$ and $n=-N,\cdots,N$,
\begin{align}
\frac{11}2y_{i}(1,n) &
-y_{i}(1,n-1)-y_{i}(1,n+1)-y_{i}(0,n)-y_{i}(2,n)\nn\\
&-\frac12y_{i}(-1,n)-y_{i}(3,n)= f_{i}, \quad i=1,2,\label{eq:inters3}
\end{align}
where $f_{1}=-\eps$ and $f_{2}=0$.

Observe
that
\[
y_{2}=0.
\]
Therefore, we only consider $y_{1}$ and omit the subscript from now on.

We first impose the Dirichlet boundary condition in the $x_1$ direction,
and the periodic boundary condition in the $x_2$ direction as
\begin{equation}\label {pbc}
\left\{
\begin{aligned}
y(-M,n)&=y(M,n)=y(M+1,n)=0,\quad n=-N,\cdots,N,\\
y(m,n)&=y(m,n+2N),\quad m=-M,\cdots,M,n=-N,\cdots,N.
\end{aligned}\right.
\end{equation}
Similar to the triangular lattice model, it is easy to check that the square lattice
model reduces to a one-dimensional chain model with the following equilibrium equations
and boundary condition.
\begin{align*}
5y(m+1)-10y(m)+5y(m-1)&=0, \; m=-M+1,\cdots,-2,\\
y(m+2)+y(m+1)-4y(m)+y(m-1)+y(m-2)&=0, \; m=2,\cdots,M-1.
\end{align*}
The equations for the interface are
\[
\left\{
\begin{aligned}
\dfrac12 y(1)+5 y(0)-\dfrac{21}{2} y(-1)+5y(-2)&=\eps,\\
y(2)+ y(1)-7 y(0)+5 y(-1)&=-2\eps,\\
y(3)+ y(2)-\dfrac72 y(1)+ y(0)+\dfrac12y(-1)&=\eps.
\end{aligned}\right.
\]
The boundary condition is
\[
y(-M)=y(M)=y(M+1)=0.
\]

It is clear that the above model is exactly the same as that has been
studied in~\cite{MingYang:2009}, which can also be obtained
from the one-dimensional model~\eqref{harmonicpotential} with
$k_1=k_2=1$.

Therefore, we impose the Dirichlet boundary condition in both
$x_1$ and $x_2$ directions as follows.
\begin{equation}\label {bc:equi}
\left\{\begin {aligned}
y(-M,n)&=y(M,n)=y(M+1,n)=0,&& n=-N,\cdots,N,\\
y(m,-N)&=y(m,N)=0,&& m=-M,\cdots,M.
\end {aligned}\right.
\end{equation}

\begin{figure}[htbp]
\centering
\subfigure[$D_{s_1}^+ y$]{%
\includegraphics[width=2.4in]{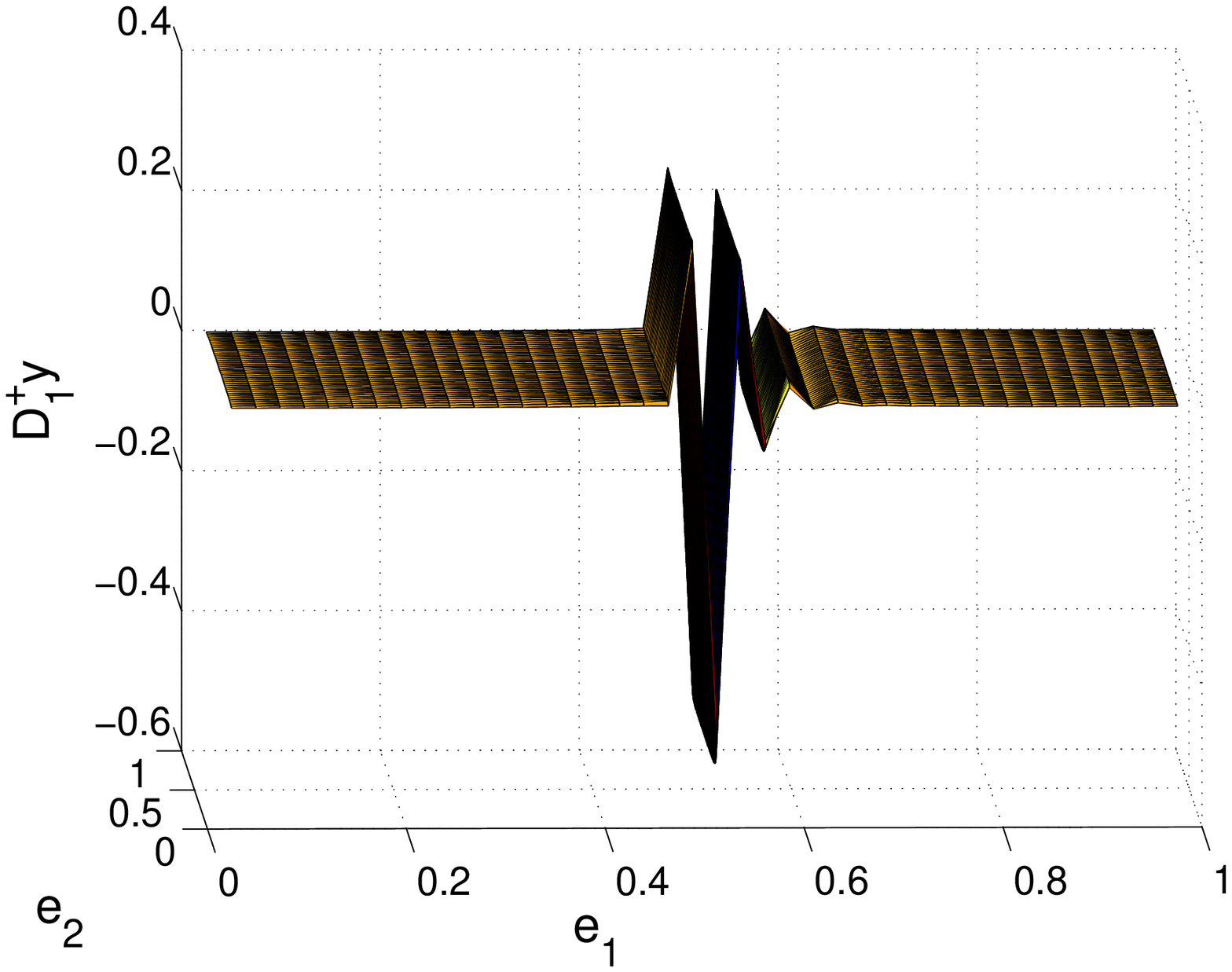}}%
\subfigure[$D_{s_2}^+ y$]{
\includegraphics[width=2.4in]{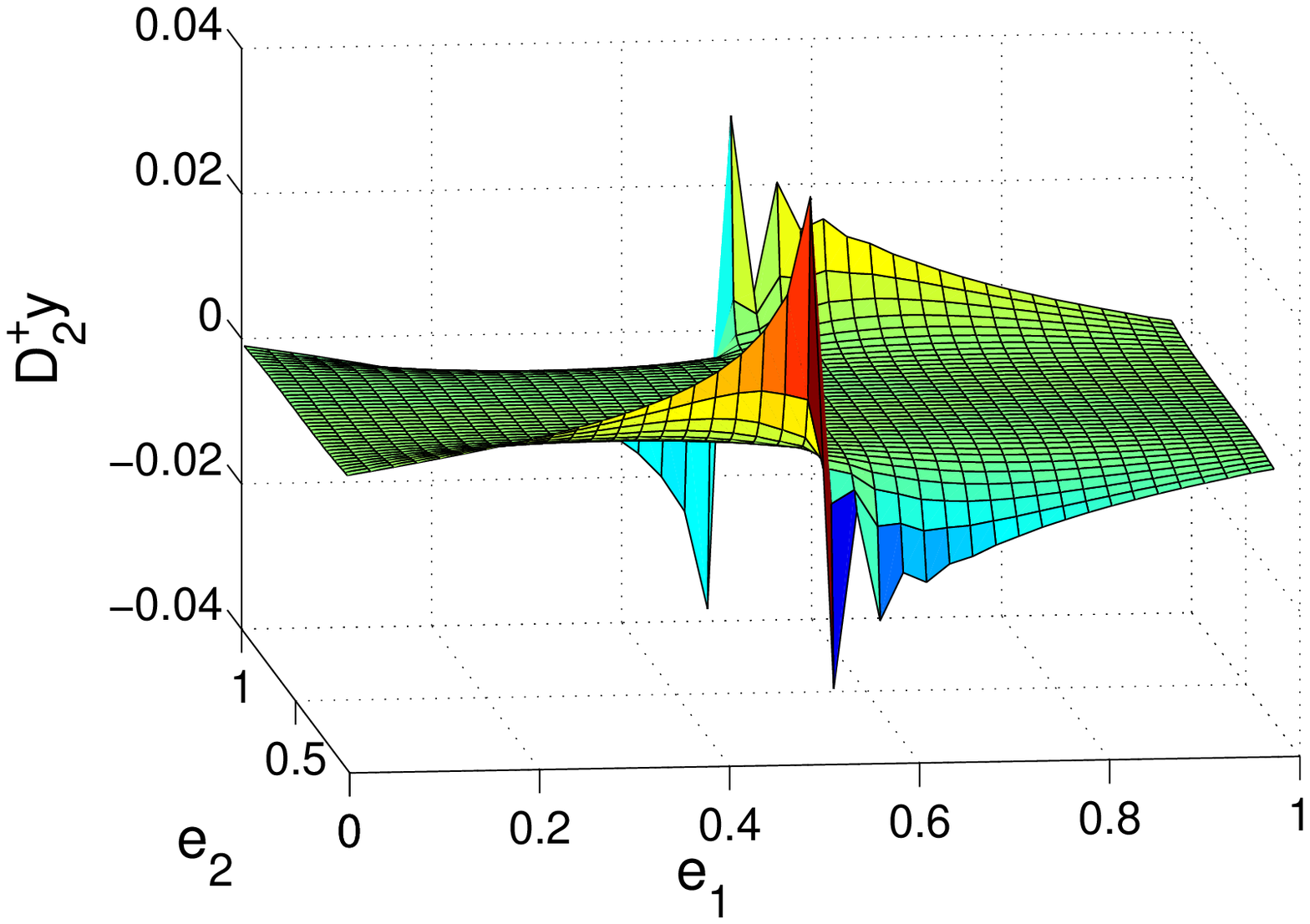}}%
\caption{\small  Profiles of the discrete gradients of $y$ for square lattice with
$M=N=20$ under Dirichlet boundary condition.}\label {Fig:der}
\end{figure}

\begin{figure}[htbp]
\centering
\subfigure[$\chi_{E_1}(x)$]{%
\includegraphics[width=2.4in]{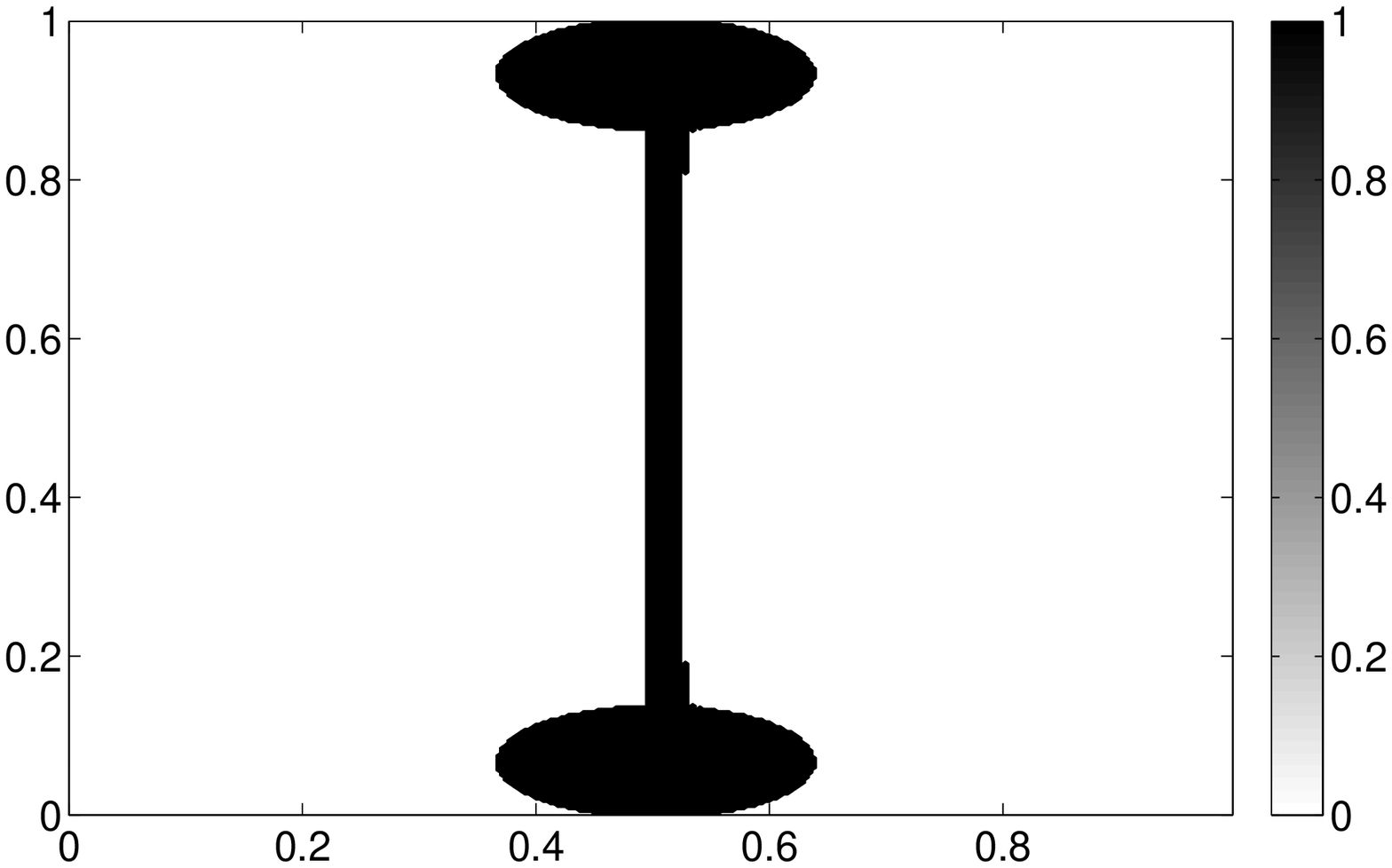}}%
\subfigure[$\chi_{E_2} (x)$]{
\includegraphics[width=2.4in]{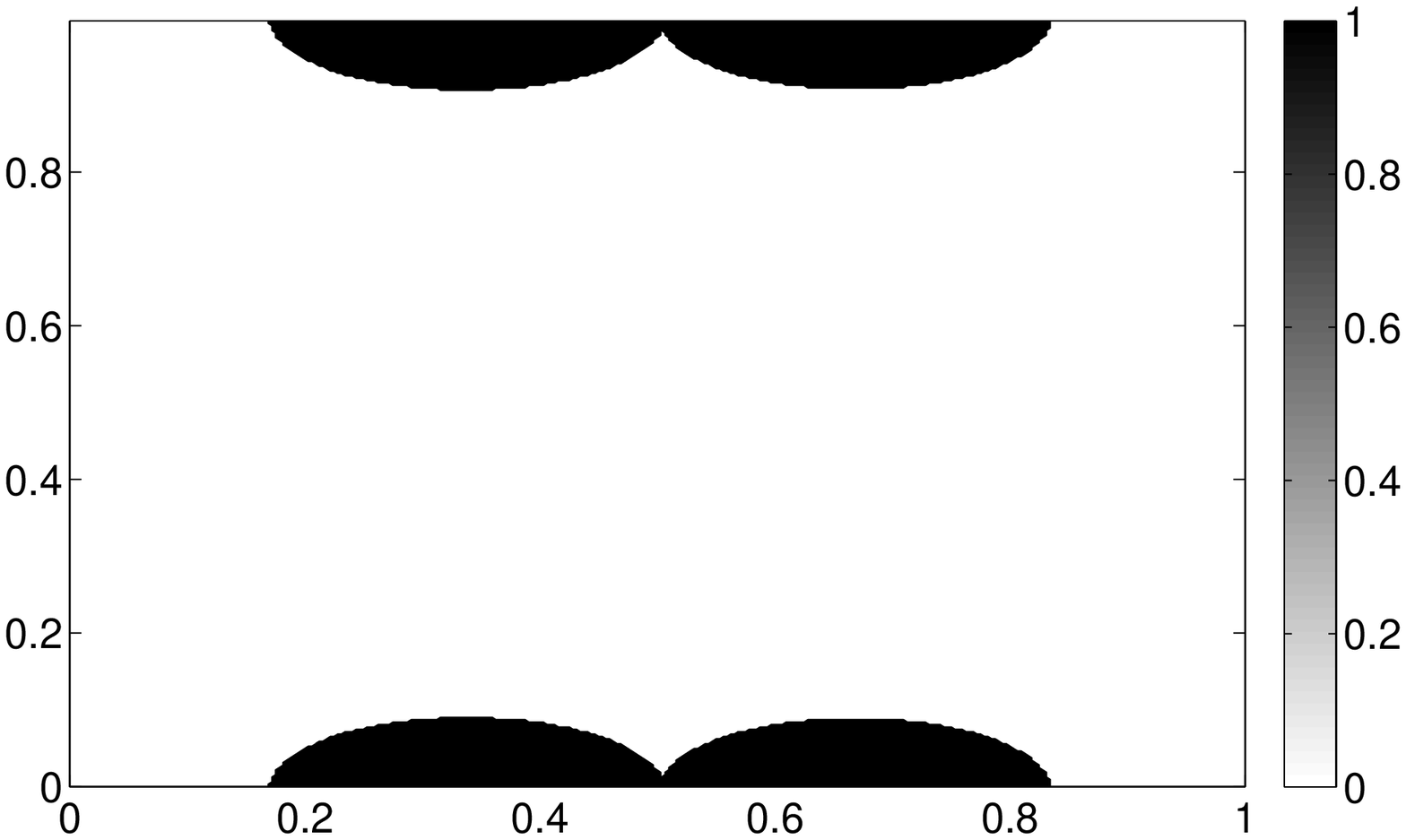}}%
\caption{\small  Profiles of the characteristic functions for square lattice.}\label {interface}
\end{figure}

Choosing $M=N=20$, we show profiles of the discrete gradients $D_{s_1}^+ y$ and $D_{s_2}^+ y$ in
Figure~\ref {Fig:der}. The feature is similar to that of the triangular lattice
model. To highlight
the interface, we plot the characteristic functions $\chi_{E_1}$ and
$\chi_{E_2}$ in Figure~\ref{interface} with $M=N=320$ and $c_0=0.04$.

We report the width of the interfacial layer in Table~\ref {layer}. It is
clear that the width of the interfacial layer is of $\co(\sqrt\eps\,)$, which is consistent with
that of the triangular lattice. In \secref{sec:estimate}, we shall prove this fact.
\begin{table}[htbp]
\centering
\begin {tabular} {ccccccc}\hline
$\eps\;(5\x10^{-2})$  & $2^0$  & $2^{-1}$ & $2^{-2}$ & $2^{-3}$ & $2^{-4}$ & \text{Rate} \\
\hline \text{Layer width}$\;(10^{-1})$  & $8$ & $5.5$ & $3.63$ & $2.56$ & $1.84$ & $0.53$\\
\hline
\end {tabular}
\vskip .3cm
\caption{\small Width of interfacial layer versus the equilibrium bond length $\eps$ for square lattice.}\label {layer}
\end {table}
\section{Exact solution for the square lattice model}
\label{sec:explicit}
To find the exact solution of the QC approximation~\eqref {eq:cb} -- \eqref{eq:inters3} with Dirichlet boundary condition~\eqref{bc:equi},
we follow the approach in~\cite {MingYang:2009}: firstly, we find
the general expression for the solution of the continuum equation
and the atomistic equation by separation of variables ansatz,
with certain unspecified constants; secondly, we use the equations
around the interface to determine these constants. The next lemma gives the general
expression of the solution.
%
\begin {lemma}\label {thm:ana}
For $m=-M,\cdots,-1$ and $n=-N,\cdots,N$, we have
\begin{equation}\label {cb-solu}
y(m,n)=\sum_{k=1}^{2N-1}a_k\sinh[(M+m)\al_k]\sin\dfrac{k\pi}{2N}(N+n),
\end{equation}
where
\[
\cosh\al_k=1+\dfrac{\lam_k}5,\qquad\lam_k=2\sin^2\dfrac{k\pi}{4N}.
\]

For $m=0,\cdots,M$ and $n=-N,\cdots,N$, we have
\begin{equation}\label {exp1}
y(m,n)=\sum_{k=1}^{2N-1}\Lr{b_k\lr{-1}^mF_m(\ga_k,\del_k)+c_kf_m(\ga_k,\del_k)}
\sin\dfrac{k\pi}{2N}(N+n),
\end{equation}
where
\begin{equation}\label {exp2}
\cosh\ga_k=\dfrac{1+\sqrt{25+8\lam_k}}4,\quad\cosh\del_k
=\dfrac{-1+\sqrt{25+8\lam_k}}4,
\end{equation}
and
\begin{equation}\label {eq:fm}
\left\{
\begin{aligned}
F_m(\ga,\del)&=\dfrac{\sinh[(M+1-m)\ga]+\sinh[(M-m)\ga]\cosh\del}
{\cosh\ga+\cosh\del}\\
&\quad-\lr{-1}^m\dfrac{\cosh[(M-m)\del]\sinh\ga}{\cosh\ga+\cosh\del},\\
f_m(\ga,\del)&=F_m(\del,\ga).
\end{aligned}\right.
\end{equation}
The coefficients $b_k$ and $c_k$ are parameters to be determined;
See; cf.,~\eqref {coef2}.
\end {lemma}

\begin {proof} By separation of variables, we get~\eqref {cb-solu}.

The explicit expression for the solution of the atomistic model can also
be obtained by separation of variables ansatz.
Substituting $y(m,n)=f(m)g(n)$ into~\eqref {eq:atom}, we get
\[
\sum_{i=-2}^2[f(m+i)-f(m)]g(n)+\sum_{i=-1}^1[g(n+i)-g(n)]f(m)=0.
\]
By~\eqref {bc:equi}, we have
\[
g(-N)=g(N)=0\quad\text{and}\quad f(M)=f(M+1)=0.
\]
We write the above equation as
\[
\dfrac{\sum_{i=-2}^2[f(m+i)-f(m)]}{f(m)}+\dfrac{\sum_{i=-1}^1[g(n+i)-g(n)]}{g(n)}=0.
\]
For $\lam\in\R$, we get
\[
\left\{
\begin {aligned}
g(n+1)+(2\lam-2)g(n)+g(n-1)&=0,\\
f(m+2)+f(m+1)-(4+2\lam)f(m)+f(m-1)+f(m-2)&=0.
\end {aligned}\right.
\]
Using the boundary condition for $g$, i.e., $g(N)=g(-N)=0$,
we have, for any $c\in\R$,
\[
g(n)=c\sin\dfrac{k\pi}{2N}(n+N)\quad\text{and}\quad
\lam=2\sin^2\dfrac{k\pi}{4N}.
\]
The characteristic equation for $f(m)$ is:
\[
t^2+t^{-2}+t+t^{-1}-2(\lam+2)=0.
\]
Denote the roots of the above
equation by $t_1,\cdots,t_4$. It is clear that
\[
t_1=-e^\ga,\;t_2=-e^{-\ga},\;t_3=e^{\del},\;t_4=e^{-\del},
\]
with
\[
\begin{aligned}
2\cosh\ga&=-s_1,&&s_1=\dfrac{-1-\sqrt{25+8\lam}}2,\\
2\cosh\del&=s_2,&&s_2=\dfrac{-1+\sqrt{25+8\lam}}2.
\end{aligned}
\]
This leads to
\begin{align*}
f(m)&=a\lr{-1}^m\sinh[(M-m)\ga]+b\lr{-1}^m\cosh[(M-m)\ga]\\
&\quad+c\sinh[(M-m)\del] +d\cosh[(M-m)\del]
\end{align*}
with constants $a,b,c$ and $d$ that
will be determined by the following conditions
\[
f(M)=f(M+1)=0.
\]
Since $M$ is even, by $f(M)=0$, we obtain
\[
b=-d.
\]
By $f(M+1)=0$, we obtain
\[
a\sinh\ga-c\sinh\del-b(\cosh\ga+\cos\del)=0.
\]
Therefore,
\[
b=\dfrac{a\sinh\ga-c\sinh\del}{\cosh\ga+\cosh\del}.
\]
We write $f(m)$ as
\begin{align*}
f(m)&=a\lr{-1}^m\sinh[(M-m)\ga]+c\sinh[(M-m)\del]\\
&\quad+\dfrac{a\sinh\ga-c\sinh\del}
{\cosh\ga+\cosh\del}\Bigl\{\lr{-1}^m\cosh[(M-m)\ga]-\cosh[(M-m)\del]\Bigr\}.
\end{align*}
It is easy to rewrite $f(m)$ into a symmetrical form
\[
f(m)=a\lr{-1}^m F_m(\ga,\del)+cF_m(\del,\ga),
\]
where $F_m(\ga,\del)$ is given in~\eqref {eq:fm}, this gives~\eqref
{exp2}.
\end {proof}

\begin{remark}
The exact solution based on the series expansion is common in finite difference.
We refer to~\cite {McCreaWhipple:1940} and~\cite {HenryBatchelor:2003}
for a thorough discussion.
\end{remark}

Next we use the interfacial equations~\eqref{eq:inters1} -- \eqref{eq:inters3}
to determine the coefficients $a_k,b_k$ and $c_k$. Denote
\begin{align*}
\yycb(m,n)&=y(m,n)\qquad-M\le m\le -1,-N\le n\le N,\\
\yat(m,n)&=y(m,n)\qquad 0\le m\le M,-N\le n\le N.
\end{align*}
Though $\yycb$ and $\yat$ are only defined piecewisely,
they actually
admit a trivial extension for any $(m,n)\in\Z^2$ as
\begin{equation}\label {observe}
\F_\eps[\yycb](m,n)=0,\qquad\fat[\yat](m,n)=0,\qquad\text{for}\quad(m,n)\in\Z^2,
\end{equation}
where
\begin{align*}
\F_\eps[y](m,n)&\equiv 12y(m,n)
-y(m,n-1)-y(m,n+1)-5y(m-1,n)-5y(m+1,n),\\
\fat[y](m,n)&\equiv 6y(m,n)
-y(m,n-1)-y(m,n+1)-y(m-1,n)-y(m+1,n)\\
&\quad-y(m-2,n)-y(m+2,n).
\end{align*}
This observation is crucial to simplify the equations around
the interface.

The equation for $m=-1$ changes to
\[
\F_\eps[\yycb(\bar1,n)]+5(\yycb-\yat)(0,n)+\frac12\Bigl(\yycb(\bar1,n)-\yat(1,n)\Bigr)
=-\eps.
\]
Using
\[
\F_\eps[\yycb(\bar1,n)]=0,
\]
we have
\begin{equation}\label {eq:inter1}
5(\yycb-\yat)(0,n)+\frac12\Bigl(\yycb(\bar1,n)-\yat(1,n)\Bigr)=-\eps.
\end{equation}

Proceeding in the same fashion, we get
\begin{align}
\yat(\bar1,n)-\frac12\yat(1,n)-\frac12\yycb(\bar1,n)=-\eps,\label {eq:inter2}\\
3\yat(0,n)+\yat(\bar1,n)+\yat(\bar2,n)-5\yycb(\bar1,n)=2\eps.\label
{eq:inter3}
\end{align}

In what follows, we use the above simplified interfacial equations and
the representation formulas to determine $a_k,b_k$ and $c_k$.

Subtracting~\eqref {eq:inter2} from~\eqref {eq:inter1}, we obtain
\[
5\yycb(0,n)+\yycb(\bar1,n)=5\yat(0,n)+\yat(\bar1,n).
\]

Substituting~\eqref {cb-solu} and~\eqref {exp1} into the above
equation, we get
\begin{equation}\label {coef1}
a_k=\dfrac{b_k(5F_0-F_{\bar1})(\ga_k,\del_k)+c_k(5f_0+f_{\bar1})(\ga_k,\del_k)}
{5\sinh[M\al_k]+\sinh[(M-1)\al_k]}.
\end{equation}

Substituting~\eqref {cb-solu} and~\eqref {exp1} into~\eqref{eq:inter2}, we get
\[
\sum_{k=1}^{2N-1}\ell_k\sin\llr{\dfrac{k\pi}{2N}(n+N)}=-2\eps
\]
with
\[
\ell_k=-\sinh[(M-1)\al_k]a_k+(-2F_{\bar1}+F_1)(\ga_k,\del_k)b_k
+(2f_{\bar1}-f_1)(\ga_k,\del_k)c_k.
\]
Using the discrete Fourier transform, we get
\[
\ell_k=\dfrac{2\x(-2\eps)}{2N-1+1}\sum_{j=1}^{2N-1}\sin\dfrac{k\pi
j}{2N}=\left\{
\begin{aligned}
-\dfrac{2\eps}N\cot{\dfrac{k\pi}{4N}},&\quad \text{if $k$ is odd, }\\
0,&\quad \text{if $k$ is even.}
\end{aligned}\right.
\]
This leads to
\begin{equation}\label {eq:bc1}
P_kb_k+p_kc_k=\ell_k,
\end{equation}
where
\begin{equation}\label{eq:expp}
\left\{
\begin {aligned}
P_k&=[-2F_{\bar1}+F_1-\rho_k(-F_{\bar1}+5F_0)](\ga_k,\del_k),\\
p_k&=[ 2f_{\bar1}-f_1-\rho_k( f_{\bar1}+5f_0)](\ga_k,\del_k),\\
\rho_k&=\dfrac{\sinh[(M-1)\al_k]}{5\sinh[M\al_k]+\sinh[(M-1)\al_k]}.
\end {aligned}\right.
\end{equation}

Using~\eqref {eq:inter2} and~\eqref {eq:inter3} to eliminate
$\yycb(\bar1,n)$, we obtain
\[
3\yat(0,n)-9\yat(\bar1,n)+\yat(\bar2,n)+5\yat(1,n)=12\eps.
\]
The coefficients $b_k$ and $c_k$ satisfy
\begin{equation}\label {eq:bc2}
R_kb_k+r_kc_k=-6\ell_k,
\end{equation}
where
\[
\left\{
\begin {aligned}
R_k&=(3F_0+9F_{\bar1}+F_{\bar2}-5F_1)(\ga_k,\del_k),\\
r_k&=(3f_0-9f_{\bar1}+f_{\bar2}+5f_1)(\ga_k,\del_k).
\end {aligned}\right.
\]
To solve the linear system~\eqref{eq:bc1} and~\eqref{eq:bc2}, we need to check
whether $P_kr_k-p_kR_k$ is nonzero for all $k$. We shall prove
in Lemma~\ref{lemma:lowdet} that this is indeed the case. Therefore, we may
solve~\eqref {eq:bc1} and~\eqref {eq:bc2} to obtain
\begin{equation}\label {coef2}
b_k=\dfrac{r_k+6p_k}{P_kr_k-p_kR_k}\ell_k,\qquad
c_k=-\dfrac{R_k+6P_k}{P_kr_k-p_kR_k}\ell_k.
\end{equation}
Substituting the above equation into~\eqref {coef1}, we get
\[
a_k=\dfrac{(r_k+6p_k)(5F_0-F_{\bar1})(\ga_k,\del_k)-(R_k+6P_k)
(5f_0+f_{\bar1})(\ga_k,\del_k)}{(5\sinh[M\al_k]+\sinh[(M-1)\al_k])
\lr{P_kr_k-p_kR_k}}\ell_k.
\]

To sum up, we obtain the representation formulas for the solution of the
QC approximation by specifying the parameters $a_k,b_k$ and $c_k$ in Lemma~\ref{thm:ana}.
\begin{theorem}\label {thm:mainexp}
Let $y$ be the solution of~\eqref{eq:cb} -- \eqref{eq:inters3}. Then for $m=-M,\cdots,\bar1$ and $n=-N,\cdots,N$,
\begin{align}
y(m,n)&=-\dfrac{2\eps}{N}\sum_{\substack{k=1\\
k\;\text{odd}}}^{2N-1} \dfrac{\Q_k}
{P_kr_k-R_kp_k}\dfrac{\sinh[(M+m)\al_k]}{\sinh[(M-1)\al_k]}\rho_k\nn\\
&\phantom{-\dfrac{2\eps}{N}\sum_{k=1}^{2N-1}}
\qquad\x\cot{\dfrac{k\pi}{4N}}\sin\llr{\dfrac{k\pi}{2N}(n+N)},\label
{cb-solu-refine}
\end{align}
where $\rho_k$ is given in~\eqref{eq:expp}$_3$.

For $m=0,\cdots,M$ and $n=-N,\cdots,N$,
\begin{equation}
y(m,n)=-\dfrac{2\eps}{N} \sum_{\substack{k=1\\
k\;\text{odd}}}^{2N-1}
\dfrac{\Q_{m,k}}{P_kr_k-R_kp_k}\cot{\dfrac{k\pi}{4N}}
\sin\llr{\dfrac{k\pi}{2N}(n+N)},\label {at-solu-refine}
\end{equation}
where
\begin{align*}
P_kr_k-R_kp_k&=6(8\rho_k-1)\begin {vmatrix}F_0&-F_{\bar1}\\
f_0&f_{\bar1}\end {vmatrix}
+(25\rho_k-3)\begin {vmatrix}-F_1&F_0\\
f_1&f_0\end {vmatrix}-\begin {vmatrix}-F_1&F_{\bar2}\\
f_1&f_{\bar2}\end {vmatrix}\\
&\quad+(2-\rho_k)\begin {vmatrix}-F_{\bar1}&F_{\bar2}\\
f_{\bar1}&f_{\bar2}\end {vmatrix}+(5\rho_k-1)\begin {vmatrix}-F_1&-F_{\bar1}\\
f_1&f_{\bar1}\end {vmatrix}-5\rho_k\begin {vmatrix}F_0&F_{\bar2}\\
f_0&f_{\bar2}\end {vmatrix},
\end{align*}
and
\begin{align*}
\Q_k&=12\begin {vmatrix} F_0 & -F_{\bar1}\\
f_0 & f_{\bar1} \end {vmatrix}
+5\begin {vmatrix}F_0&F_{\bar2}+F_1\\
f_0&f_{\bar2}-f_1\end {vmatrix}
+\begin {vmatrix}-F_{\bar1}&F_{\bar2}+F_1\\
f_{\bar1}&f_{\bar2}-f_1\end {vmatrix},
\end{align*}
and
\begin{align*}
\Q_{m,k}
&=3(1-10\rho_k)\begin {vmatrix}\lr{-1}^m F_m&F_0\\
f_m& f_0\end {vmatrix}\\
&\quad+3(1-2\rho_k)\begin {vmatrix}\lr{-1}^m F_m& -F_{\bar1}\\
f_m&f_{\bar1}\end {vmatrix}+\begin {vmatrix}\lr{-1}^m F_m&F_{\bar2}+F_1\\
f_m&f_{\bar2}-f_1\end {vmatrix}.
\end{align*}
\end{theorem}

As an immediate consequence of the above theorem, the solution is symmetrical
with respect to $n=0$, i.e.,
\begin{equation}\label{eq:sym}
y(m,n)=y(m,-n),
\end{equation}
which can be easily verified from the representation
formulas~\eqref{cb-solu-refine} and~\eqref{at-solu-refine}.
\section{Estimate of the QC Solution for the square lattice model}
\label{sec:estimate}
The main result of this section is the following pointwise estimate of the solution.
\begin {theorem}\label {main}
Let $y$ be the solution of~\eqref{eq:cb} -- \eqref{eq:inters3}. Then
\begin{align}
\abs{Dy(m,n)}&\le\dfrac{C}{m^2}\exp\left[-\dfrac{\abs{m}}{6\sqrt{5}N}\right],\;
m\le-1,\label {cb-solu-final}\\
\abs{D y(m,n)}&\le
C\Lr{\Lr{\dfrac{3-\sqrt{5}}2}^m+\dfrac{1}{m^2+1}\exp\left[-\dfrac{2m}{15N}\right]},\; m\ge
0.\label {at-solu-final}
\end{align}
\end {theorem}

A direct consequence of the above theorem is the estimate of the width of the interfacial
layer, that is, the region beyond which $\abs{D y}=\co(\eps)$.
\begin {corollary}\label {main2}
The width of the interfacial layer is $\co(\sqrt{\eps})$.
\end {corollary}

We exploit the explicit expression of the solution in Theorem~\ref {thm:mainexp} to prove
Theorem~\ref {main}. Note that the
terms $P_kr_k-R_kp_k,\Q_k$ and $\Q_{m,k}$ consist of
the terms like $(-1)^mF_mf_n-(-1)^nF_nf_m$ for
different integers $m$ and $n$. The asymptotical behavior
of such terms will be given in a series of lemmas, i.e., Lemma~\ref {lemma:lowdet}, Lemma~\ref{lemma:lead1} and Lemma~\ref{lemma:lead2}.
We begin with certain elementary estimates that will be frequently used later on.
\begin {lemma}
For $1\le k\le 2N-1$, there holds
\begin{align}
\dfrac{\lam_k}6&\le\cosh\ga_k-\dfrac32\le\dfrac{\lam_k}5,\label {est-gamma}\\
\sinh\del_k&\ge\sqrt{\dfrac{\lam_k}3},\label {eq:estdel}\\
\sinh\al_k&\ge\sqrt{\dfrac{2\lam_k}5},\quad\sinh\dfrac{\al_k}2
=\sqrt{\dfrac{\lam_k}{10}}.\label {eq:estal}
\end{align}
\end {lemma}

\begin {proof}
Invoking~\eqref {exp2}, we have
\[
\cosh\ga_k-\dfrac32=\dfrac14\Lr{\sqrt{25+8\lam_k}-5}
=\dfrac{2\lam_k}{\sqrt{25+8\lam_k}+5},
\]
which immediately implies~\eqref {est-gamma}.

The estimate~\eqref{eq:estal} follows from~\eqref {cb-solu} by definition.

Using~\eqref{exp2}, we have
\begin{equation}\label{eq:gadel}
\cosh\gamma_k-\cosh\del_k=\dfrac12,
\end{equation}
which together with the definition leads to
\begin{align*}
\sinh^2\del_k&=\cosh^2\del_k-1=(\cosh\gamma_k-1/2)^2-1\\
&=(\cosh\gamma_k+1/2)(\cosh\gamma_k-3/2)\\
&=(\cosh\gamma_k-3/2+2)(\cosh\gamma_k-3/2).
\end{align*}
Using~\eqref{est-gamma}, we have
\[
\sinh^2\del_k\ge 2(\cosh\gamma_k-3/2)\ge\dfrac{\lam_k}3.
\]
This gives~\eqref {eq:estdel}.
\end {proof}

To proceed further, we need the following estimates.
\begin {lemma}
For $1\le k\le 2N-1$, there holds
\begin{align}
\exp(-\al_k)&\le\exp\Bigl[-\dfrac{k}{2\sqrt{5}N}\Bigr],\label
{exp-est2}\\
\exp(-\ga_k)&\le\dfrac{3-\sqrt{5}}2,\quad
\exp(-\del_k)\le\exp\Bigl[-\dfrac{k}{5N}\Bigr].\label
{exp-est1}
\end{align}
\end {lemma}

\begin {proof}
We only prove~\eqref {exp-est2}. Other cases are similar.

Using~\eqref {eq:estal} and $\cosh\al_k\ge 1$, we have
\begin{align*}
\exp\al_k&=\cosh\al_k+\sinh\al_k\\
&\ge 1+\sqrt{\dfrac{2\lam_k}5}\\
&=1+\dfrac{2}{\sqrt5}\sin\dfrac{k\pi}{4N}.
\end{align*}
Using Jordan's inequality
\[
\sin x\ge\dfrac{2}{\pi}x\quad\text{for\quad}x\in[0,\pi/2],
\]
we have
\[
\exp\al_k\ge
1+\dfrac{k}{\sqrt{5}N}.
\]
For any $0<x<2/\sqrt{5}$, we have
\[
\ln(1+x)\ge x(1-x/2)\ge x(1-1/\sqrt{5})\ge x/2.
\]
Using the fact that $k/(\sqrt{5}N)\le 2/\sqrt{5}$ since $1\le k\le 2N-1$, and
combining the above two inequalities, we obtain
\begin{align*}
\exp(-\al_k)&\le\Lr{1+\dfrac{k}{\sqrt{5}N}}^{-1}
=\exp[-\ln(1+k/(\sqrt{5}N))]\\
&\le\exp\left[-\dfrac{k}{2\sqrt{5}N}\right].
\end{align*}
\end {proof}

The next lemma concerns the estimate of $\rho_k$.
\begin {lemma}
\begin{equation}\label {est-rho}
0<1-6\rho_k\le\dfrac56\Lr{\dfrac{\lam_k}5+\dfrac{1}{M-1}+\sinh\al_k}.
\end{equation}
\end {lemma}

\begin{proof}
Using the definition of $\rho_k$, we get
\[
1-6\rho_k=\dfrac{5\Lr{\sinh[M\al_k]-\sinh[(M-1)\al_k]}}
{5\sinh[M\al_k]+\sinh[(M-1)\al_k]},
\]
which implies the left hand side of~\eqref {est-rho}. Moreover
\begin{align*}
1-6\rho_k&\le\dfrac56\Lr{\dfrac{\sinh[M\al_k]}{\sinh[(M-1)\al_k]}-1}\\
&=\dfrac56\Lr{\cosh\al_k-1+\cot[(M-1)\al_k]\sinh\al_k}. \end{align*}
Using $\cosh t\le 1+\sinh t$ for any $t\in\R$, we have
\begin{align*}
\cot[(M-1)\al_k]\sinh\al_k&\le\sinh\al_k+\dfrac{\sinh\al_k}{\sinh[(M-1)\al_k]}\\
&\le\sinh\al_k+\dfrac{1}{M-1},
\end{align*}
where we have used the elementary inequality
\[
\dfrac{\sinh[Mt]}{\sinh t}\ge M.
\]
Combining the above three inequalities, we obtain the right hand side of~\eqref {est-rho}.
\end {proof}

By the definition of $\rho_k$ and the left hand side of~\eqref {est-rho}, we get
\begin{equation}\label{eq:rho}
0<\rho_k\le 1/6.
\end{equation}

A direct calculation gives
\begin{align*}
(\cosh\ga+\cosh\del)&\bigl((-1)^mF_mf_n-(-1)^nF_nf_m\bigr)\\
&=A\sinh[M\ga]\sinh[M\del]+B\cosh[M\ga]\sinh[M\del]\nn\\
&\quad+C\sinh[M\ga]\cosh[M\del]+D\cosh[M\ga]\cosh[M\del]\nn\\
&\quad-\sinh[(m-n)\del]\sinh\ga+(-1)^{m+n}\sinh[(m-n)\ga]\sinh\del,
\end{align*}
where
\begin{align*}
A&=(-1)^m\cosh[m\ga]\cosh[(n-1)\del]+(-1)^m\cosh[n\del]\cosh[(m-1)\ga]\\
&\quad-(-1)^n\cosh[n\ga]\cosh[(m-1)\del]-(-1)^n\cosh[m\del]\cosh[(n-1)\ga],
\end{align*}
and
\begin{align*}
B&=-(-1)^m\sinh[m\ga]\cosh[(n-1)\del]-(-1)^m\sinh[(m-1)\ga]\cosh[n\del]\\
&\quad+(-1)^n\sinh[n\ga]\cosh[(m-1)\del]+(-1)^n\sinh[(n-1)\ga]\cosh[m\del],
\end{align*}
and
\begin{align*}
C&=-(-1)^m\cosh[m\ga]\sinh[(n-1)\del]-(-1)^m\cosh[(m-1)\ga]\sinh[n\del]\\
&\quad+(-1)^n\cosh[n\ga]\sinh[(m-1)\del]+(-1)^n\cosh[(n-1)\ga]\sinh[m\del],
\end{align*}
and
\begin{align*}
D&=(-1)^m\sinh[m\ga]\sinh[(n-1)\del]+(-1)^m\sinh[(m-1)\ga]\sinh[n\del]\\
&\quad-(-1)^n\sinh[n\ga]\sinh[(m-1)\del]-(-1)^n\sinh[(n-1)\ga]\sinh[m\del].
\end{align*}

The following lemma gives a lower bound for $\abs{P_kr_k-R_kp_k}$.
\begin{lemma}\label{lemma:lowdet}
There holds
\begin{equation}\label {low-det}
(\cosh\ga_k+\cosh\del_k)\abs{P_kr_k-R_kp_k}\ge
\dfrac{5}{24}\Lr{1-\exp\left[-\dfrac{2M}{5N}\right]}\exp{[M(\ga_k+\del_k)]}.
\end{equation}
\end{lemma}

\begin {proof}
A direct calculation gives
\begin {align*}
&\quad(\cosh\ga_k+\cosh\del_k)(P_kr_k-R_kp_k)\\
&=A_k\sinh[M\ga_k]
\sinh[M\del_k]+B_k\cosh[M\ga_k]\sinh[M\del_k]\\
&\quad+C_k\sinh[M\ga_k]\cosh[M\del_k]+D_k\cosh[M\ga_k]\cosh[M\del_k]\\
&\quad+2(12+2\lam_k-72\rho_k)\sinh\ga_k\sinh\del_k,
\end {align*}
where
\begin {align*}
A_k&=\Lr{4\lam_k^2+\dfrac{313}2\lam_k+315}\rho_k
           -6\lam_k^2-39\lam_k-\dfrac{105}2,\\
B_k&=\Lr{18\lam_k+\dfrac{225}4+\Lr{\dfrac{229}4+2\lam_k}\sqrt{25+8\lam_k}}
           \rho_k\sinh\ga_k\\
&\quad+\Lr{\dfrac{25}2+4\lam_k-(14+4\lam_k)\sqrt{25+8\lam_k}}\sinh\ga_k,\\
C_k&=\Lr{-18\lam_k-\dfrac{225}4+\Lr{\dfrac{229}4+2\lam_k}\sqrt{25+8\lam_k}}
           \rho_k\sinh\del_k\\
&\quad+\Lr{-\dfrac{25}2-4\lam_k-(14+4\lam_k)\sqrt{25+8\lam_k}}\sinh\del_k,\\
D_k&=\Lr{(169+8\lam_k)\rho_k-(74+20\lam_k)}\sinh\ga_k\sinh\del_k.
\end {align*}

Using~\eqref{eq:rho}, we may show
\begin{equation}\label{eq:neg}
A_k, B_k, C_k, D_k < 0.
\end{equation}
Using $D_k<0$ and invoking~\eqref {eq:rho} once again, we have
\begin{align}
&\quad D_k\cosh[M\ga_k]\cosh[M\del_k]
+2(12+2\lam_k-72\rho_k)\sinh\ga_k\sinh\del_k\nn\\
&\le D_k+2(12+2\lam_k-72\rho_k)\sinh\ga_k\sinh\del_k\nn\\
&=\Lr{(169+8\lam_k)\rho_k-(74+20\lam_k)+2(12+2\lam_k-72\rho_k)}
\sinh\ga_k\sinh\del_k\nn\\
&=\Lr{(24+8\lam_k)\rho_k-(50+16\lam_k)}
\sinh\ga_k\sinh\del_k\nn\\
&\le -(46+44\lam_k/3)\sinh\ga_k\sinh\del_k\nn\\
&<0,\label{eq:combneg}
\end{align}
which together with~\eqref{eq:neg} implies
\[
(\cosh\ga_k+\cosh\del_k)(P_kr_k-R_kp_k)<0.
\]
Combining the above equation with~\eqref{eq:neg} and~\eqref{eq:combneg},
we obtain
\begin{align}
&\labs{(\cosh\ga_k+\cosh\del_k)(P_kr_k-R_kp_k)}
=-(\cosh\ga_k+\cosh\del_k)(P_kr_k-R_kp_k)\nn\\
&=-B_k\cosh[M\ga_k]\sinh[M\del_k]\nn\\
&\quad-A_k\sinh[M\ga_k]\sinh[M\del_k]-C_k\sinh[M\ga_k]\cosh[M\del_k]\nn\\
&\quad+\Lr{-D_k\cosh[M\ga_k]\cosh[M\del_k]
-2(12+2\lam_k-72\rho_k)\sinh\ga_k\sinh\del_k}\nn\\
&\ge\abs{B_k}\cosh[M\ga_k]\sinh[M\del_k].\label{eq:lowmed}
\end{align}

Using~\eqref{est-rho}, we bound $\abs{B_k}$ as
\begin{align*}
\abs{B_k}&=-B_k\nn\\
&=-\Lr{18\lam_k+\dfrac{225}4+\Lr{\dfrac{229}4+2\lam_k}\sqrt{25+8\lam_k}}
           \rho_k\sinh\ga_k\nn\\
&\quad-\Lr{\dfrac{25}2+4\lam_k-(14+4\lam_k)\sqrt{25+8\lam_k}}\sinh\ga_k\nn\\
&\ge\Lr{\Lr{\dfrac{107}{24}+\dfrac{11}{3}\lam_k}
\sqrt{25+8\lam_k}-\dfrac{175}8-7\lam_k}\sinh\ga_k.
\end{align*}
A direct calculation gives $\sinh\ga_k\ge\sqrt5/2$, which together with
the above inequality yields
\begin{equation}\label{eq:lowb}
\abs{B_k}\ge\dfrac56.
\end{equation}

By~\eqref{exp-est1},
\[
\exp[2M\del_k]\ge\exp[2kM/(5N)]\ge\exp[2M/(5N)].
\]
It follows from the above inequality that
\begin{align*}
\sinh[M\del_k]&=\dfrac12\exp[M\del_k]\Lr{1-\exp[-2M\del_k]}\\
&\ge\dfrac12\Lr{1-\exp[-2M/(5N)]}\exp[M\del_k].
\end{align*}
Substituting the above inequality and~\eqref{eq:lowb}
into~\eqref{eq:lowmed} implies~\eqref {low-det}.
\end {proof}

Next two lemmas concern the upper bounds of $\Q_k$ and $\Q_{m,k}$. Instead
of calculating all the coefficients of $\Q_k$ and $\Q_{m,k}$ as we have done
for $(\cosh\ga_k+\cosh\del_k)\abs{P_kr_k-R_kp_k}$, we consider the coefficients
of the leading order terms of $\Q_k$ and $\Q_{m,k}$. We write
\begin{align}\label {sepa}
(\cosh\ga+\cosh\del)&\bigl((-1)^mF_mf_n-(-1)^nF_nf_m\bigr)\nn\\
&=\dfrac14(A+B+C+D)e^{M(\ga+\del)}\nn\\
&+\dfrac14(-A-B+C+D)e^{M(\ga-\del)}+\lot
\end{align}
with
\begin{equation}\label {4sum}
\left\{
\begin{aligned}
A+B+C+D&=(e^{\ga}+e^{\del})\Lr{(-1)^me^{-(n\del+m\ga)}-(-1)^ne^{-(n\ga+m\del)}},\\
-A-B+C+D&=(e^{\ga}+e^{-\del})\Lr{-(-1)^me^{n\del-m\ga}+(-1)^ne^{m\del-n\ga}},
\end{aligned}\right.
\end{equation}
and $\lot$ stands for the terms that are of lower order than $e^{M(\ga+\del)}$
and $e^{M(\ga-\del)}$.

Using~\eqref {cb-solu-refine}, we write
$(\cosh\ga_k+\cosh\del_k)\Q_k$ as
\[
(\cosh\ga_k+\cosh\del_k)\Q_k=\Q_k^0\exp[M(\ga_k+\del_k)]
+\Q_k^1\exp[M(\ga_k-\del_k)]+\lot
\]
with
\begin{align*}
\Q_k^0&=\dfrac14(e^{\ga_k}+e^{\del_k})\Bigl\{12(e^{\ga_k}+e^{\del_k})
+5(e^{2\del_k}-e^{2\ga_k})
-5(e^{-\ga_k}+e^{-\del_k})\\
&\phantom{(e^{\ga_k}+e^{\del_k})\Bigl\{12}
\quad-e^{\ga_k+\del_k}(e^{\ga_k}+e^{\del_k})+e^{\ga_k-\del_k}-e^{\del_k-\ga_k}\Bigr\},\\
\Q_k^1&=\dfrac14(e^{\ga_k}+e^{-\del_k})\Bigl\{-12(e^{\ga_k}+e^{-\del_k})
+5(e^{2\ga_k}-e^{-2\del_k})
+5(e^{-\ga_k}+e^{\del_k})\\
&\phantom{(e^{\ga_k}+e^{-\del_k})\Bigl\{12}
\qquad+e^{\ga_k-\del_k}(e^{\ga_k}+e^{-\del_k})
+e^{-\ga_k-\del_k}-e^{\ga_k+\del_k}\Bigr\}.
\end{align*}
We have the following estimate for $\Q_k^0$.
\begin{lemma}\label{lemma:lead1}
There exists $C$ such that
\begin{equation}\label {lead-1}
\abs{\Q_k^0}+\abs{\Q_k^1}\le C\sqrt{\lam_k}.
\end{equation}
\end {lemma}

\begin {proof}
We only estimate $\Q_k^0$, and $\Q_k^1$ can be bounded
similarly. We firstly write $\Q_k^0$ as
\begin{align}
\Q_k^0&=\dfrac14(e^{\ga_k}+e^{\del_k})\Bigl\{(12+e^{-\del_k}-e^{2\del_k})e^{\ga_k}
-(5+e^{\del_k})(e^{2\ga_k}+e^{-\ga_k})\nn\\
&\phantom{(e^{\ga_k}+e^{\del_k})\Bigl\{12}
\quad+12e^{\del_k}-5e^{-\del_k}+5e^{2\del_k}\Bigr\}.\label {qk}
\end{align}
By definition, we have
\begin{equation}\label {rela3}
\cosh[2\ga_k]=2\cosh^2\ga_k-1=\cosh\ga_k+\lam_k+2,
\end{equation}
which implies
\begin{align*}
e^{2\ga_k}+e^{-\ga_k}&=\cosh[2\ga_k]+\sinh[2\ga_k]+\cosh\ga_k-\sinh\ga_k\\
&=2\cosh\ga_k+\lam_k+2+\sinh\ga_k(2\cosh\ga_k-1)\\
&=2(e^{\ga_k}+1)+\lam_k+\sinh\ga_k(2\cosh\ga_k-3).
\end{align*}
Substituting the above equation into~\eqref {qk} produces
\begin{align*}
\Q_k^0&=\dfrac14(e^{\ga_k}+e^{\del_k})
\Bigl\{(5-e^{\ga_k})[2(e^{\del_k}-1)-e^{-\del_k}+e^{2\del_k})]\\
&\phantom{(e^{\ga_k}+e^{\del_k})\Bigl\{}
\qquad-(5+e^{\del_k})[\lam_k+\sinh\ga_k(2\cosh\ga_k-3)]\Bigr\}.
\end{align*}
Using~\eqref{eq:gadel}, we get
\begin{align*}
2(e^{\del_k}-1)-e^{-\del_k}+e^{2\del_k}&=(e^{\del_k}-1)(3+2\cosh\del_k)\\
&=2(\cosh\ga_k+1)(e^{\del_k}-1).
\end{align*}
Combining the above two equations, we obtain
\begin{align}
\Q_k^0&=(e^{\ga_k}+e^{\del_k})\Bigl\{2(5-e^{\ga_k})(\cosh\ga_k+1)(e^{\del_k}-1)\nn\\
&\phantom{(e^{\ga_k}+e^{\del_k})\Bigl\{}\qquad
-(5+e^{\del_k})[\lam_k+\sinh\ga_k(2\cosh\ga_k-3)]\Bigr\}.\label {qk1}
\end{align}

Proceeding along the same way that leads to~\eqref {qk1}, we obtain
\begin{align}
\Q_k^1&=\dfrac14(e^{\ga_k}+e^{\del_k})\Bigl\{
(5+e^{-\del_k})[\lam_k+2\sinh\ga_k(2\cosh\ga_k-3)]\nn\\
&\phantom{\dfrac14(e^{\ga_k}+e^{\del_k})}
\qquad+(5-e^{\ga_k})(1+3e^{-\del_k}+e^{-2\del_k})(e^{\del_k}-1)\Bigr\}.\label {qk2}
\end{align}

Using~\eqref{eq:gadel} gives
\[
e^{\del_k}-1=\cosh\del_k+\sinh\del_k-1=\cosh\ga_k-\dfrac32+\sinh\del_k.
\]

Substituting the above equation into~\eqref {qk1} and~\eqref {qk2}, and
using the estimates~\eqref{est-gamma} and~\eqref {eq:estdel}, we
obtain~\eqref {lead-1}.
\end {proof}

Next we write $(\cosh\ga_k+\cosh\del_k)\Q_{m,k}$ as
\begin{align*}
(\cosh\ga_k+\cosh\del_k)\Q_{m,k}&=\Q_{m,k}^0\exp[M(\ga_k+\del_k)]\\
&\quad+\Q_{m,k}^1\exp[M(\ga_k-\del_k)]+\text{L.O.T.}
\end{align*}
with
\[
\left\{
\begin {aligned}
\Q_{m,k}^0&=\dfrac14(e^{\ga_k}+e^{\del_k})\Lr{(-1)^me^{-m\ga_k}\Q_1-e^{-m\del_k}\Q_2},\\
\Q_1&=3(1+e^{\del_k})+e^{2\del_k}-e^{-\del_k}-6\rho_k(5+e^{\del_k}),\\
\Q_2&=e^{2\ga_k}+e^{-\ga_k}-3(e^{\ga_k}-1)+6\rho_k(e^{\ga_k}-5),
\end {aligned}
\right.
\]
and
\[
\left\{
\begin {aligned}
\Q_{m,k}^1&=\dfrac14(e^{\ga_k}+e^{-\del_k})\Lr{-(-1)^me^{-m\ga_k}\Q_3-e^{m\del_k}\Q_4},\\
\Q_3&=3(1-2\rho)(e^{-\del_k}-1)+e^{-2\del_k}-e^{\del_k}+6(1-6\rho_k),\\
\Q_4&=-\Q_2.
\end {aligned}
\right.
\]

We have the following estimate for $\Q_{m,k}^0$ and $\Q_{m,k}^1$ .
\begin{lemma}\label{lemma:lead2}
There holds
\begin{equation}\label {lead-2}
\abs{\Q_1}+\abs{\Q_2}+\abs{\Q_3}+\abs{\Q_4}\le C\Lr{\sqrt{\lam_k}+\dfrac{1}{M-1}}.
\end{equation}
\end{lemma}

\begin {proof}
We only estimate $\Q_1$ and $\Q_2$. The terms $\Q_3$ and
$\Q_4$ can be bounded similarly.

Similar to~\eqref {rela3}, we have
\[
\cosh[2\del_k]=-\cosh\del_k+\lam_k+2.
\]
Using the above equation, we write $\Q_1$ as
\begin{align*}
\Q_1&=3(1+\cosh\del_k+\sinh\del_k)-\cosh\del_k+\lam_k+2
+\sinh[2\del_k]\\
&\quad-(\cosh\del_k-\sinh\del_k)-6\rho_k(5+e^{\del_k})\\
&=\lam_k+5+\cosh\del_k+4\sinh\del_k+\sinh[2\del_k]-6\rho_k(5+e^{\del_k})\\
&=\lam_k+(5+e^{\del_k})(1-6\rho_k)+\sinh\del_k(3+2\cosh\del_k).
\end{align*}

Using~\eqref {rela3} and proceeding along the same line that leads to the above expression of
$\Q_1$, we obtain
\begin{align*}
\Q_2&=2\cosh^2\ga_k+2+\sinh[2\ga_k]-3(\cosh\ga_k+\sinh\ga_k)+\cosh\ga_k-\sinh\ga_k\\
&\quad+6\rho_k(e^{\ga_k}-5)\\
&=\lam_k+5-\cosh\ga_k+\sinh[2\ga_k]-4\sinh\ga_k+6\rho_k(e^{\ga_k}-5)\\
&=\lam_k+(1-6\rho_k)(5-e^{\ga_k})+2\sinh\ga_k(\cosh\ga_k-3/2).
\end{align*}
Using~\eqref {est-gamma},~\eqref {eq:estdel} and~\eqref {est-rho}, we get
\[
\abs{\Q_1}+\abs{\Q_2}\le C\Lr{\sqrt{\lam_k}+\dfrac{1}{M-1}}.
\]
\end {proof}

To prove Theorem~\ref {main}, we need the following identity that can be found in~\cite [p. 38, formula 1.353(1)]{JeffreyZwillinger:2007}.
\begin {lemma}\label {tri-iden}
For any $\vr\in(0,1)$, we have
\[
\sum_{k=1}^N\vr^{2k-1}\sin[(2k-1)x]=\dfrac{(\vr+\vr^3)\sin
x-\vr^{2N+1}\sin[(2N+1)x]+\vr^{2N+3}\sin[(2N-1)x]}{1-2\vr^2\cos[2x]+\vr^4}.
\]
\end {lemma}

Based on the above estimates, we are ready to prove Theorem~\ref {main}.
\vskip .5cm
\noindent {\em Proof of Theorem~\ref {main}\;}
Using~\eqref {low-det} and~\eqref {lead-1}, we have, for $m\le -2$,
\begin{equation}
\abs{Dy(m,n)}\le\dfrac{C}{N}\sum_{\substack{k=1\\
k\;\text{odd}}}^{2N-1}\exp(-\abs{m}\al_k)\sin\dfrac{k\pi}{2N}.\label
{cb-solu-det}
\end{equation}
By~\eqref {low-det} and~\eqref {lead-2}, we have, for $m\ge 0$,
\begin{align}
\abs{D_1y(m,n)}\le\dfrac{C}{N}\sum_{\substack{k=1\\
k\;\text{odd}}}^{2N-1}\Lr{\exp(-m\ga_k)+\exp(-m\del_k)\sin\dfrac{k\pi}{2N}},\label
{at-solu-det1}\\
\abs{D_2y(m,n)}\le\dfrac{C}{N}\sum_{\substack{k=1\\
k\;\text{odd}}}^{2N-1}\Lr{\exp(-m\ga_k)\sin\dfrac{k\pi}{2N}
+\exp(-m\del_k)\sin\dfrac{k\pi}{2N}}.\label {at-solu-det2}
\end{align}

Let $\vr=\exp[-\abs{m}/(2\sqrt{5}N)]$ and $x=\pi/[2N]$, and using
Lemma~\ref {tri-iden}, we have
\begin{align*}
\sum_{\substack{k=1\\
k\;\text{odd}}}^{2N-1}\exp(-\abs{m}\al_k)
\sin\dfrac{k\pi}{2N}&\le\sum_{k=1}^N\vr^{2k-1}\sin[(2k-1)x]\\
&=\dfrac{\vr(1+\vr^2)(1+\vr^{2N})}{(1-\vr^2)^2+4\vr^2\sin^2x}\sin x\\
&\le\dfrac{2\vr(1+\vr)^2}{(1-\vr^2)^2}\sin x\\
&=\dfrac{2\vr}{(1-\vr)^2}\sin x.
\end{align*}
Using Lozarevi\'c's inequality~\cite {Lazarevic:1966}:
\[
\cosh t\le\Lr{\dfrac{\sinh t}{t}}^3, \qquad t\not=0,
\]
and the elementary inequality
$\cosh t\ge e^t/2,t\in\R$,
we obtain
\begin{align*}
\dfrac{2\vr}{(1-\vr)^2}\sin
x&=\dfrac{\sin\dfrac{\pi}{2N}}{2\sinh^2\Lr{\dfrac{\abs{m}}{4\sqrt{5}N}}}\\
&\le\dfrac{\dfrac{\pi}{2N}}
{\Lr{\dfrac{\abs{m}}{4\sqrt{5}N}}^2\cosh^{2/3}\Lr{\dfrac{\abs{m}}{4\sqrt{5}N}}}\\[.1in]
&\le\dfrac{2^{2/3}40\pi N}{m^2}\exp\left[-\dfrac{\abs{m}}{6\sqrt{5}N}\right]\\[.1in]
&\le\dfrac{80\pi
N}{m^2}\exp\left[-\dfrac{\abs{m}}{6\sqrt{5}N}\right],
\end{align*}
which together with~\eqref {cb-solu-det} leads to\eqref
{cb-solu-final}.

For $m\ge 1$ and let $\vr=\exp[-m/(5N)]$, we immediately have~\eqref
{at-solu-final}.

The proof of the case when $m=-1$ can be done in the same
way that leads to~\eqref {cb-solu-final}. We leave it to the
interested readers. \qed

Proceeding along the same line that leads to Theorem~\ref {main}, we have the
following estimate on the solution.
\begin {corollary}
 There exists $C$ such that
 \begin{align*}
   \abs{y(m,n)}&\le C\dfrac{\eps}{\abs{m}},\; m\le -1,\\
   \abs{y(m,n)}&\le C\eps\Lr{\Lr{\dfrac{3-\sqrt5}2}^m+\dfrac{1}{m+1}},\; m\ge 0.
   \end{align*}
 \end {corollary}
The above estimate suggests that the ghost force actually induces a negligible error on
the solution, which is as small as $\eps$.
\bibliographystyle{amsplain}
\bibliography{plan}
\end {document}